\documentclass[preprint,3p,11pt]{elsarticle}

\usepackage{color}
\usepackage{amssymb}
\usepackage{amsmath}
\usepackage{amsthm}

\usepackage{epsfig}
\usepackage{ulem}
\usepackage{epic}
\usepackage{tabls}
\usepackage{booktabs}
\usepackage{threeparttable}
\usepackage{pict2e}

\usepackage{float}
 \usepackage{subfig}

\usepackage{hyperref}

\hypersetup{backref,colorlinks=true,linkcolor=blue}

\usepackage[vlined]{algorithm2e}

\begin{document}
\begin{frontmatter}



\title{
 Nonlinear Iterative Projection  Methods with  Multigrid in Photon Frequency  for Thermal Radiative Transfer
}

\author{Dmitriy Y. Anistratov}
\address{Department of Nuclear Engineering,
North Carolina State University, Raleigh, NC 27695\\
anistratov@ncsu.edu}

\begin{abstract}

This paper presents  nonlinear  iterative methods
for the fundamental  thermal radiative transfer (TRT) model defined by the time-dependent multifrequency radiative transfer (RT) equation
and the material energy balance (MEB) equation.
The iterative methods are based on the nonlinear projection approach and use  multiple grids  in photon frequency.
 They are formulated by  the high-order RT equation on a given grid in photon frequency
 and   low-order moment equations on a hierarchy of frequency grids.
 The material temperature is evaluated in the subspace of the lowest dimensionality
  from the   MEB equation  coupled to the effective grey low-order equations.
 The   algorithms apply various multigrid cycles to visit  frequency grids.
Numerical results are presented to demonstrate convergence of the multigrid iterative algorithms in
TRT problems with large number of photon frequency groups.
 \end{abstract}


\end{frontmatter}

\section{Introduction}

Radiation transport is of fundamental  importance for  high-temperature phenomena in
high-energy density physics, inertial confinement fusion, astrophysics etc \cite{shu-1991,drake-2006,Graziani-2004}.
High-energy photons interact with matter  and play  essential role in redistribution of energy in a physical system.
The dynamics of the system is affected by absorption and emission of photons by matter
  leading to change in  the material temperature.
The propagation of photons and their interaction with matter is described by the Boltzmann  radiative transfer (RT) equation in which the opacities and emission source term are highly nonlinear functions of
  the material temperature \cite{Chandr-1960}.
The basic thermal radiative transfer (TRT)  model is defined by
the RT equation and  the material energy balance (MEB) equation.
The material energy depends nonlinearly on temperature.
The advanced model is formulated by the system of radiation hydrodynamics equations
that includes conservation of mass and momentum equations  accounting for the effects of radiation
\cite{zel-1966,mihalas-FRH-1984}.

The dimensionality of the TRT  problem
is driven by the specific intensity of radiation that depends on
time, spatial position, direction of particle motion, and photon frequency.
In a general geometry, it is a seven-dimensional function.
The parameters of matter depend only on space and time and hence four independent variables.
The TRT system of equations is approximated
 by implicit methods to avoid stability constraints on the   time step size.
 This yields a complicated and large system of equations on  the phase-space grid.
At a time level, the RT equation  can be solved deterministically
 by  sweeping  a spatial mesh in each discrete direction  for a photon frequency group.
 The  costs of transport calculations depend  on
   complexity of a transport discretization scheme  and
  the degree of  phase space resolution by  the grid.
The  computational effort on the same spatial grid  scales with increase in number of angular directions and frequency groups.
The RT equation is tightly coupled  to the MEB equation.
  Simple iterations  between the RT and MEB equations converge very slow.
It is necessary to apply efficient iterative techniques to reduce number of transport iterations that involve transport sweeps on  the phase-space grid.

To solve the nonlinear TRT problem,  the   RT and MEB equations can be  linearized with respect to temperature.
This yields  a modified high-order  RT equation with pseudo-scattering in the  phase space on every time step \cite{fleck-1971,larsen-jcp-1988}.
To accelerate transport iterations  associated  with   pseudo-scattering,
  fast  iteration methods, such as, the diffusion-synthetic acceleration, grey transport acceleration, Krylov methods
are  used \cite{larsen-jcp-1988,Alcouffe-85,morel-jqsrt-1985,morel-jcp-2007,mla-m&c2013}.

The nonlinear projection approach (NPA)  is based on  applying  projection operators in angular and frequency variables
 to formulate a system of low-order equations for moments of the intensity \cite{gol'din-cmmp-1964,Goldin-sbornik-82}.
 The low-order moment equations are closed with the high-order RT equation
  by means of  prolongation operators and exact closures.
 The iteration methods for solving RT equation based on the NPA are
 the quasidiffusion (QD) method  (aka Eddington Variable Factor (VEF) method),
 nonlinear diffusion acceleration, coarse-mesh finite differencing (CMFD) method, partial current-based  CMFD,
 $\alpha$-weighted methods \cite{gol'din-cmmp-1964,auer-mihalas-1970,dya-ewl-2001,lr-dya-2010,kord-nda-1984,kord-physor-2002,cho-lee-jkns-2003}.
 The NPA has been applied to TRT and radiation hydrodynamics problems
 \cite{gol'din-1972,surzh-82,PASE-1986,Winkler-85,dya-aristova-vya-mm1996,aristova-vya-avk-m&c1999,park-et-al-2012,yee-2016,dya-jcp-2019}.

A group of  iterative algorithms for frequency-dependent radiative transfer
use  two  frequency grids involving the grey  low-order transport problem
 which is coupled  to the MEB equation to estimate the material temperature during iterations
 \cite{gol'din-1972,dya-aristova-vya-mm1996,aristova-vya-avk-m&c1999,dya-jcp-2019}.
 The grey low-order problem is derived by a nonlinear projection
of   multigroup moment  equations in photon frequency
using  the group radiation energy density and fluxes to define grey coefficients.
The computational methods based on the  linearized high-order RT and MEB equations
formulate a grey problem for the modified RT equation with pseudo-scattering
using the grey opacities averaged with slowest converging  error modes
\cite{larsen-jcp-1988,morel-jqsrt-1985,morel-jcp-2007}.

In this study, the TRT  problems with  large number of groups are considered.
It has been demonstrated that iterative methods based on nonlinear projection  and  multiple grids in particle energy
are efficient algorithms for particle transport problems
\cite{lrc-dya-nse-2016,lrc-dya-ks-nse-2019,lc-dya-mc2019,lee-2012}.
In this paper, we present   new nonlinear iterative methods for TRT problems with multigrid
in photon frequency.  They are formulated by  means of the high-order RT equation
 on a given grid in photon frequency  and   low-order QD   equations on multiple   frequency grids.
 The hierarchy  of grids for the low-order equations always include (i) the given fine grid as the first one and
(ii)  the grid with one interval that covers the whole frequency range.
The MEB equation is coupled to the low-order equations on the coarsest frequency grid, namely, to the grey low-order equations.  As a result, the material temperature is evaluated in the subspace of the lowest dimensionality.

The reminder of the paper is organized as follows.
In Sec. \ref{sec:trt}, the TRT problem is described.
In Sec. \ref{sec:npm},  the multilevel QD (MLQD) method with two frequency grids is  reviewed.
In Sec. \ref{sec:mlqd-mgrid}, the  MLQD method with multiple grids in photon frequency is   formulated in continuous form; different multigrid algorithms are   presented.
The discretization of equations of the   MLQD method is formulated in  Sec. \ref{sec:disc}.
The numerical results are presented in  Sec \ref{sec:res}.
In Sec. \ref{sec:concl}, we conclude with a brief discussion.

\section{\label{sec:trt}Thermal Radiative Transfer Problem}

We consider the TRT  model in one-dimensional slab geometry.
It  is defined by the frequency-dependent RT equation \cite{zel-1966}
\begin{equation} \label{rt-nu}
\frac{1}{c} \frac{\partial I}{\partial t}(x, \mu, \nu,  t)
+ \mu \frac{\partial I}{\partial x}(x, \mu, \nu,  t)
+ \sigma(\nu, T)  I(x, \mu, \nu,  t)
= \sigma(\nu, T)   B(\nu,T) \, ,
\end{equation}
\begin{equation*}
x \in [0, X] \, , \quad    \mu \in [-1, 1] \, , \quad   t \ge 0 \, , \quad     \nu \in [0, \infty) \, ,
\end{equation*}
and  the MEB equation
\begin{equation}\label{eb-nu}
 \frac{\partial  \varepsilon(T)}{\partial t} =
 \int_0^{\infty}\int_{-1}^1   \sigma(\nu, T)  \Big( I(x, \mu, \nu,  t) -    B(\nu,T) \Big) d \mu   d \nu \,
\end{equation}
with the initial conditions
\begin{equation} \label{rt-nu-ic}
 I|_{t=0}=  I^0 \, ,
\end{equation}
\begin{equation}
T|_{t=0}= T^0  \, ,
\end{equation}
and the boundary conditions
\begin{subequations} \label{rt-nu-bc}
\begin{equation}
 I|_{x=0}=  I^+\, , \ \mu \in (0,1]  \, ,
\end{equation}
\begin{equation}
  I|_{x=X}=  I^- \, , \ \mu \in [-1,0)  \, .
\end{equation}
\end{subequations}
Here
 $I_{\nu}$ is the  specific intensity,
 $T$ is the material temperature;
$\varepsilon$ is the material energy density;
$\sigma_{\nu}$ is the photon opacity;
$x$ is the spatial position;
$\mu$ is     the directional cosine of particle motion;
$\nu$ is the photon frequency;
$t$ is time.
\begin{equation}
B(\nu, T)=\frac{4\pi h\nu^3}{c^2}\frac{1}{e^{\frac{h\nu}{kT}}-1}
\end{equation}
 is the  Planck black-body distribution function multiplied by $2\pi$,
 where $h$ is the Planck’s constant, $c$ is the speed of light, $k$ is the Boltzmann’s constant.
 The  TRT model \eqref{rt-nu} and   \eqref{eb-nu}  neglects material motion, scattering, heat conduction, and external sources. It is applicable  in the case of the supersonic radiation wave \cite{moore-2015}.

To formulate discretization of   the  RT equation with respect to the frequency variable,
we  define  the grid
\begin{equation} \label{fine-grid}
\Omega_{\nu}=\{\nu_g,  g \in \mathbb{N}(n_{\nu}+1) \} \,
\end{equation}
that divides the whole frequency range  into discrete groups given by intervals
$\omega_g=[\nu_g, \nu_{g+1}]$.
Here $n_{\nu}$ is the number of groups, $\mathbb{N}(n_{\nu}+1)=\{1,\ldots,n_{\nu}+1\}$,
$\nu_1=0$, and $\nu_{n_{\nu}+1} =\infty $ or some  maximum value.
The TRT  model in the multigroup approximation is defined by the   RT equation
\begin{equation} \label{rt}
\frac{1}{c} \frac{\partial I_g}{\partial t} + \mathcal{L}_g I_g =  q_g \, , \quad
 g \in \mathbb{N}(n_{\nu}) \, ,
\end{equation}
 \begin{equation}\label{rt-bc}
\left. I_g \right|_{x=0}  = I_g^{+} \, , \   \mu > 0 \, ,
\quad
\left. I_g \right|_{x=X}  = I_g^{-} \, , \   \mu < 0 \, ,
\end{equation}
\begin{equation}\label{rt-ic}
\left. I_g \right|_{t=0 }
 = I_g^0 \,
\end{equation}
 for the group intensity
 \begin{equation}
 I_g = \int_{\nu_g}^{\nu_{g+1}} I d \nu \, ,
 \end{equation}
where
\begin{equation}
 \mathcal{L}_g I_g \equiv \mu \frac{\partial I_g}{\partial x} + \sigma_{E,g} I_g \, ,
\end{equation}
\begin{equation}
 q_g = \sigma_{B,g} B_g \, ,
 \end{equation}
\begin{equation}
B_g=\int_{\nu_g}^{\nu_{g+1}}  B d \nu \, .
\end{equation}
The    RT equation  \eqref{rt}  for the group $g$ is formulated by means of two different  group opacities:
\begin{equation}
\sigma_{B,g}(T) = \frac{\displaystyle \int_{\nu_g}^{\nu_{g+1}}   \sigma(\nu, T) B(\nu, T) d \nu}
{\displaystyle \int_{\nu_g}^{\nu_{g+1}}   B(\nu, T) d \nu} \, ,
\quad
\sigma_{E,g}(T,T_r) = \frac{\displaystyle \int_{\nu_g}^{\nu_{g+1}}   \sigma(\nu, T) B(\nu, T_r) d \nu}
{\displaystyle \int_{\nu_g}^{\nu_{g+1}}   B(\nu, T_r) d \nu} \,
\end{equation}
 averaged  with the Planck spectrum function at the material temperature $T$  and
 the effective  temperature of radiation, $T_r$, respectively
\cite{PASE-1986,dya-aristova-vya-mm1996,aristova-vya-avk-m&c1999}.
The  multigroup form of the MEB equation is given by
\begin{equation}\label{eb}
\frac{\partial \varepsilon}{\partial t} =
 \sum_{g=1}^{n_{\nu}} \int_{-1}^1 \Bigl( \sigma_{E,g}  I_g -  \sigma_{B,g}B_{g}\Big) d \mu   \, .
\end{equation}

\section{\label{sec:npm} Nonlinear Projection Methods}

 \subsection{The Two-Level QD  Method}

 To solve the TRT problem \eqref{rt} and \eqref{eb},
we apply the nonlinear projection approach
based on the QD (VEF) method \cite{gol'din-cmmp-1964,auer-mihalas-1970}.
The two-level QD  method on the given  grid in photon frequency  is formulated by means of
projection in the angular variable and exact closures of the moment equations.
The system of equations of this method consists of two  parts:
(i) the high-order RT equation  and
(ii) the multigroup low-order QD (LOQD) equations on the grid $\Omega_{\nu}$
for the angular moments of the group intensity.
 The LOQD problem for  the  group radiation energy density
  \begin{equation} \label{E_g}
E_g(x,t) = \frac{1}{c} \int_{-1}^1 I_g (x,\mu,t)d \mu\,
\end{equation}
and  flux
\begin{equation} \label{F_g}
F_g(x,t) =  \int_{-1}^1 \mu I_g(x,\mu,t) d \mu \,
\end{equation}
is defined by the  moment equations
 \begin{subequations} \label{mloqd}
\begin{equation}\label{mloqd-1}
\frac{\partial E_g   }{\partial t}
+ \frac{\partial F_g   }{\partial x}
+ c\sigma_{E,g}  E_g
= 2  \sigma_{B, g}B_g \, ,
\quad
\end{equation}
\begin{equation}\label{mloqd-2}
\frac{1}{c}\frac{\partial F_g   }{\partial t}
+ c\frac{\partial (f_g E_g)   }{\partial x}
+ \sigma_{R,g}  F_g = 0 \,    ,
\end{equation}
\end{subequations}
where
\begin{equation} \label{f_g}
f_g = \frac{\displaystyle\int_{-1}^1 \mu^2 I_g d \mu}{ \displaystyle \int_{-1}^1 I_g d \mu} \,
\end{equation}
is  the group  QD (Eddington)  factor that defines the exact closure for the LOQD equations.
The boundary and initial conditions are given by \cite{gol'din-cmmp-1964,gol'din-1972}
\begin{equation} \label{mloqd-bcs}
F_g\big|_{x=0}  = \big( c \, C_g^- (E_g - E_g^{in+}) +   F_g^{in+} \big) \big|_{x=0}\, ,
\quad
F_g\big|_{x=X} = \big( c  \, C_g^+ (E_g - E_g^{in-}) +   F_g^{in-} \big) \big|_{x=X} \, ,
\end{equation}
\begin{equation}\label{mloqd-ic}
E_g \big|_{t=0 }  = E_g^0 \, , \quad  F_g \big|_{t=0 }  = F_g^0 \, ,
\end{equation}
where
\begin{equation} \label{c_g}
C_g^-  =  \frac{\displaystyle \int_{-1}^0 \mu I_g(0,\mu,t) d \mu}{\displaystyle \int_{-1}^0 I_g (0,\mu,t)  d \mu} \, , \quad
C_g^+  = \frac{\displaystyle \int_0^1 \mu I_g(X,\mu,t) d \mu}{\displaystyle \int_0^1  I_g(X,\mu,t)  d \mu}\Bigg|_{x=X} \,
\end{equation}
are the boundary QD factors that define the closure for the energy density and flux at the boundary of the spatial domain, and
\begin{equation}
E_g^{\pm}  = \pm \frac{1}{c}\int_{0}^{\pm 1} I_g^{\pm} d \mu \, , \quad
F_g^{\pm}  = \pm \int_{0}^{\pm 1} \mu I_g^{\pm} d \mu \, ,
\end{equation}
\begin{equation}
E_g^0  = \frac{1}{c}\int_{-1}^1 I_g^0 d \mu \, , \quad
F_g^0  =  \int_{-1}^1 \mu I_g^0 d \mu \, .
\end{equation}
The first moment equation (\ref{mloqd-2}) is defined with
 the Rosseland group opacity
\begin{equation}
\sigma_{R,g}(T,T_r) =  \frac{\displaystyle \int_{\nu_g}^{\nu_{g+1}}   \frac{\partial B(\nu, T')}{\partial T'}\Big|_{T'=T_r} d \nu}
{\displaystyle \int_{\nu_g}^{\nu_{g+1}}   \frac{1}{\sigma(\nu, T)} \frac{\partial B(\nu, T')}{\partial T'}\Big|_{T'=T_r}  d \nu}
 \, .
\end{equation}
The  LOQD equations can be written in the following general operator form:
\begin{equation} \label{group-loqd}
 \frac{\partial \mathbf{Y}_g}{\partial t} + \mathcal{M}_g \mathbf{Y}_g  = \mathbf{Q}_g \, , \
\mathbf{Y}_g  =( E_g, F_{g})^T \, , \ g \in \mathbb{N}(n_{\nu}) \, ,
\end{equation}
where the operator $\mathcal{M}_g = \mathcal{M}_g[f_g,C_g^{\pm},T]$,  $f_g=f_g[I_g]$,
 $C_g^{\pm}=C_g^{\pm}[I_g]$, and $\mathbf{Q}_g=\mathbf{Q}_g[T]$.
To couple  the MEB equation with the LOQD equations, it is cast  in the multigroup form
in terms of the group energy densities as follows:
\begin{equation} \label{eb-group}
\frac{\partial \varepsilon(T)}{\partial t} =
 \sum_{g=1}^{n_{\nu}} \Big( c \sigma_{E,g}E_{g} -  2 \sigma_{B,g}B_{g} \Big) \, .
\end{equation}
In summary, the system of equations of the two-level  QD method  on the  frequency grid $\Omega_{\nu}$ is defined by Eqs. \eqref{rt}, \eqref{mloqd}, and \eqref{eb-group}.
The main feature of this method is that the estimation of temperature is performed in the projected space  by solving  the MEB equation coupled to the multigroup LOQD  equations.

Algorithm \ref{MLQD-algorithm-1} describes  elements of the iterative scheme for the two-level  QD method
at every time step. Here $s$ is the  index of transport iterations; $j$ is the index of the time step.
The first stage of the iteration algorithm is to update group opacities $\sigma_{E,g}$, $\sigma_{B,g}$
using the latest estimation of the temperature $T^{(s)}$.
This defines the operator $\mathcal{L}_g^{(s)}=\mathcal{L}_g[T^{(s)}]$.
On the second stage, the  multigroup high-order RT equations are solved to obtain
  the group intensities $I_g^{(s)}$.
Then,  $I_g^{(s)}$  is used   to compute the QD factors $f_g[I_g^{(s)}]$  and $C_g^{\pm}[I_g^{(s)}]$.
This defines the low-order operator $\mathcal{M}_g^{(s)}=\mathcal{M}_g[f_g^{(s)},C_g^{\pm(s)}, T^{(s)}]$.
 On the next stage,  the new estimation of temperature is calculated by  solving  the MEB equation coupled to the multigroup LOQD  equations defined by  the operator $\mathcal{M}_g^{(s)}$.
Various numerical techniques can be applied to solve the nonlinear system of the multigroup LOQD and MEB equations
on the frequency grid $\Omega_{\nu}$.

\vspace{0.5cm}
\begin{algorithm}[H]
\DontPrintSemicolon
$s=0$, $T^{(1)}=T^{j-1}$\;
  \While{$|| T^{(s)} - T^{(s-1)} ||>\epsilon||T^{(s)}||   \, \& \,
 ||E^{(s)} - E^{(s-1)}||> \epsilon||E^{(s)}|| $ }{
$\bullet$ transport  iterations\;
 $s=s+1$ \;
$T^{(s)}   \Rightarrow  \mathcal{L}_g^{(s)}$  \;
$c^{-1}  \partial_t I_g  + \mathcal{L}_g^{(s)} I_g = q_g[T^{(s)}]   \ \mbox{on} \  \Omega_{\nu}
\Rightarrow$  $I_g^{(s)}, g \in \mathbb{N}(n_{\nu})$ \;
$I_g^{(s)} \Rightarrow    f_{g}^{(s)}, C_g^{\pm (s)}, g \in \mathbb{N}(n_{\nu})$\;
 $T^{(s)}, f_{g}^{(s)}, C_g^{\pm (s)}  \Rightarrow  \mathcal{M}_g^{(s)}, g \in \mathbb{N}(n_{\nu})$\;
$\partial_t \mathbf{Y}_g \! +  \!  \mathcal{M}_g^{(s)} \mathbf{Y}_g \!  =  \!  \mathbf{Q}_g$ on  $\Omega_{\nu}$   \&   MEB Eq.~\eqref{eb-group}
$\Rightarrow$ $T^{(s+1)}$,  $E_g^{(s+1)}$,  $F_g^{(s+1)}$,  $g \in \mathbb{N}(n_{\nu})$\;
}
\caption{ \label{MLQD-algorithm-1}  The two-level QD method on  the single frequency grid $\Omega_{\nu}$.}
\end{algorithm}
\vspace{1cm}

\subsection{\label{sec:mlqd-2grid} The Multilevel QD  Method with Two Grids in Frequency}

In this section, we review the multilevel QD (MLQD) method that uses  two frequency grids  to solve the multigroup LOQD equations coupled with the MEB equation
\cite{PASE-1986,dya-aristova-vya-mm1996,aristova-vya-avk-m&c1999,dya-jcp-2019}.
This method introduces  the coarse frequency grid
\begin{equation} \label{grey-grid}
\Omega_{\nu}^{\ast} =\{ \nu_1=0,  \nu_2=\infty \}
\end{equation}
 with one  group.
The effective  grey  LOQD equations are formulated on the grid $\Omega_{\nu}^{\ast}$
 for the  total radiation energy density
\begin{equation}
E(x,t)= \sum_{g=1}^{n_{\nu}} E_g(x,t) \,
\end{equation}
 and total flux
 \begin{equation}
F(x,t)=  \sum_{g=1}^{n_{\nu}} F_g(x,t) \, .
\end{equation}
The effective grey LOQD problem   is defined by
\begin{subequations}\label{gloqd}
\begin{equation}\label{gloqd-1}
\frac{\partial E}{\partial t}
+ \frac{\partial F}{\partial x}
+ c \bar \sigma_E E
= c  \bar \sigma_B a_R T^4\, ,
\end{equation}
\begin{equation}\label{gloqd-2}
\frac{1}{c}\frac{\partial F }{\partial t}
+ c\frac{\partial (\bar f E)   }{\partial x}
+ \bar \sigma_R  F + \bar \eta  E  = 0 \,
\end{equation}
\end{subequations}
with the boundary conditions
\begin{equation} \label{gloqd-bcs}
F\big|_{x=0}  = \big( c \, \bar C^- (E - E^{+}) +   F^{+} \big) \big|_{x=0}\, ,
\quad
F\big|_{x=X} = \big( c  \, \bar C^+ (E - E^{-}) +   F^{-} \big) \big|_{x=X} \,
\end{equation}
and the initial conditions
\begin{equation}\label{gloqd-ic}
E \big|_{t=0 }  = E^0 \, , \quad  F \big|_{t=0 }  = F^0 \, ,
 \end{equation}
 where
 \begin{equation}
E^{\pm}= \sum_{g=1}^{n_{\nu}} E_g^{\pm} \, ,
\quad
F^{\pm}= \sum_{g=1}^{n_{\nu}} F_g^{\pm} \,  ,
\quad
E^0= \sum_{g=1}^{n_{\nu}} E_g^0 \, ,
\quad
F^0= \sum_{g=1}^{n_{\nu}} F_g^0 \, ,
\end{equation}
 $a_R$ is the Stefan's constant.
 The grey LOQD equations \eqref{gloqd} on the coarse  grid  $\Omega_{\nu}^{\ast}$ are derived
 by   summing the group LOQD equations \eqref{mloqd} over all  groups
and  formulating  exact closures by means of   the grey QD factor
\begin{equation}
\bar f = \big<f \big>_{E}  \, ,
\end{equation}
the grey  boundary QD factors
\begin{equation}
\bar  C^- = \big< C^- \big>_{E}\big|_{x=0}  \, ,
\quad
\bar  C^+ = \big< C^+ \big>_{E}\big|_{x=X}  \, ,
\end{equation}
the  grey opacities
 \begin{equation}
 \bar \sigma_E = \big< \sigma_{E}\big>_{E}  \, , \quad
  \bar \sigma_B = \big< \sigma_{B}\big>_{B}  \, , \quad
   \bar \sigma_R = \big< \sigma_{R}\big>_{|F|}  \, , \quad
 \end{equation}
where the notations for the averaged qualities  are  defined by
 \begin{equation}
 \big<\psi \big>_{H}   = \frac{\sum_{g=1}^{n_{\nu}} \psi_g H_g}
 {\sum_{g=1}^{n_{\nu}} H_g} \, , \quad
H_g = \begin{cases}
E_g \ \mbox{for} \   H=E  \, , \\
B_g \ \mbox{for} \ H=B  \, ,  \\
|F_g| \ \mbox{for}  \ H=|F|    \, .
\end{cases}
\end{equation}
The compensation term is given by
  \begin{equation}
\bar \eta  = \frac{\sum_{g=1}^{n_{\nu}} (\sigma_{R,g} - \bar \sigma_R) F_g}
{\sum_{g=1}^{n_{\nu}} E_g} \, .
\end{equation}
The operator  form of Eqs. \eqref{gloqd} is the following:
\begin{equation} \label{grey-loqd}
 \frac{\partial \mathbf{Y}}{\partial t} +  \mathcal{\bar M} \mathbf{Y} = \mathbf{Q} \, ,
\quad
\mathbf{Y} = ( E, F)^T \, ,
\end{equation}
where the operator $\mathcal{\bar  M}  = \mathcal{\bar M}[f_g,C_g^{\pm},E_g,F_g,T]$ and $\mathbf{Q}=\mathbf{Q}[T]$.
The effective grey LOQD equations
are coupled with the   MEB  equation in the grey form
\begin{equation}\label{eb-grey}
\frac{\partial \varepsilon(T)}{\partial t} =  c  \big(\bar \sigma_E E - \bar\sigma_B a_R T^4\big) \, .\
\end{equation}
In summary, the system of equations of the MLQD method  on  two  frequency grids $\Omega_{\nu}$ and  $\Omega_{\nu}^{\ast}$
 is defined by Eqs. \eqref{rt}, \eqref{mloqd}, \eqref{gloqd}, and \eqref{eb-grey}.

Algorithm \ref{mlqd-2grids} shows the iterative scheme for the two-grid  MLQD method
on the $j$-th time step.
 The iterative scheme  consists of  nested iterations.
 Here $\ell$ is the  index of inner (low-order) iterations.
The outer iteration cycle is the transport iteration.
Note that there are no transport sweeps for $s=0$.
On each outer  iteration, the system of the group  LOQD equations
on the given  grid in frequency $\Omega_{\nu}$  and   MEB equation
is solved by means of the grey LOQD equations.
The estimation of temperature is obtained from   Eqs. \eqref{gloqd} and \eqref{eb-grey}.

\vspace{0.5cm}
\begin{algorithm}[H]
\DontPrintSemicolon
$s=-1$, $T^{(0)}=T^{j-1}$, $f_{g}^{(0)}=f_{g}^{j-1}$\;
\While{$|| T^{(s)} - T^{(s-1)} ||>\epsilon||T^{(s)}|| \, \& \,
||E^{(s)} - E^{(s-1)}||> \epsilon||E^{(s)}|| $ }{
$\bullet$ transport (outer) iterations\;
 $s=s+1$ \;
 \If{$s>0$}{
$T^{(s)} \Rightarrow \mathcal{L}_g^{(s)}$ \;
$c^{-1}  \partial_t I_g  + \mathcal{L}_g^{(s)} I_g = q_g[T^{(s)}]   \ \mbox{on} \  \Omega_{\nu}
\Rightarrow I_g^{(s)}, g \in \mathbb{N}(n_{\nu})$\;
 $ I_g^{(s)} \Rightarrow f_{g}^{(s)}, C_g^{\pm (s)}, g \in \mathbb{N}(n_{\nu})$  \;
 }
$\ell=-1$,  $T^{(1,s)}=T^{(s)} $ \;
\While{$|| T^{(\ell,s)} - T^{(\ell-1,s)} ||>\tilde\epsilon||T^{(\ell,s)}||  \, \& \,
 || E^{(\ell,s)} - E^{(\ell-1,s)} ||>\tilde\epsilon||E^{(\ell,s)}||$ }{
 $\ell=\ell+1$ \;
$ T^{(\ell,s)}, f_g^{(s)},C_g^{\pm (s)} \Rightarrow \mathcal{M}_g^{(\ell,s)}$\;
$\partial _t \mathbf{Y}_g+ \mathcal{M}_g ^{(\ell,s)}\mathbf{Y}_g = \mathbf{Q}_g[T^{(\ell,s)}]    \  \mbox{on} \ \Omega_{\nu}\Rightarrow  E_g^{(\ell,s)}, F_g^{(\ell,s)}, g \in \mathbb{N}(n_{\nu}) $ \;
$T^{(\ell,s)}, E_g^{(\ell,s)}, F_g^{(\ell,s)}, f_g^{(s)},C_g^{\pm (s)} \Rightarrow \mathcal{\bar M}^{(\ell,s)}$\;
$\partial _t \mathbf{Y}+ \mathcal{\bar M}^{(\ell,s)} \mathbf{Y} = \mathbf{Q}$ on $\Omega_{\nu}^{\ast}$
\&  MEB Eq. \eqref{eb-grey} $\Rightarrow T^{(\ell+1,s)}, E^{(\ell+1,s)}, F^{(\ell+1,s)}$\;
}
$T^{(s+1)} \leftarrow  T^{(\ell+1,s)}$ \;
}
\caption{\label{mlqd-2grids}  The MLQD method with two  frequency grids  $\Omega_{\nu}$ and $\Omega_{\nu}^{\ast}$.}
\end{algorithm}

\section{\label{sec:mlqd-mgrid} The Multilevel QD Method with Multigrid in Frequency}

We  now define a sequence  of nested grids in frequency
\begin{equation}
\mathcal{G}_{\Omega}^{\Gamma} = \{  \Omega_{\nu}^{\gamma} , \gamma \in \mathbb{N}(\Gamma) \} \, , \quad
 \Omega_{\nu}^{\gamma}=\{\nu_p^{\gamma},  p \in \mathbb{N}(n_{\nu}^{\gamma}+1)  \} , \quad
  \Omega^{1}_{\nu}  \equiv  \Omega_{\nu} \, ,  \quad
   \Omega^{\Gamma}_{\nu}    \equiv      \Omega_{\nu}^{\ast} \, ,
\end{equation}
where $\gamma$ is the grid index, and $p$ is the index of the frequency interval $\omega^{\gamma}_p=[ \nu^{\gamma}_p,  \nu^{\gamma}_{p+1}]$.
These grids are defined by successive coarsening. The $p$-th interval
 of the grid  $\Omega_{\nu}^{\gamma}$
is formed by  intervals of the grid $\Omega_{\nu}^{\gamma-1}$ with the set of indices
 $\Lambda_p^{\gamma} = \{ p': \omega^{\gamma-1}_{p'} \in \omega^{\gamma}_p \}$
and  hence
$\omega^{\gamma}_p = \underset{ p' \in \Lambda_p^{\gamma}}{\bigcup} \omega^{\gamma-1}_{p '}=\underset{ g \in \tilde \Lambda_p^{\gamma}}{\bigcup} \omega_{g}^1$,
where  $\tilde  \Lambda_p^{\gamma} = \{ g: \omega^{1}_{g} \in \omega^{\gamma}_p \}$.
The group LOQD equations on the coarse  grid $\Omega_{\nu}^{\gamma}$ for
\begin{equation}
E_p^{\gamma} =  \int_{\nu^{\gamma}_p}^{\nu^{\gamma}_{p+1}}  E_{\nu}d \nu  \, ,
\quad
F_p^{\gamma} =\int_{\nu^{\gamma}_p}^{\nu^{\gamma}_{p+1}}  F_{\nu}d \nu \,
\end{equation}
are derived by projecting the group LOQD equations from the fine grid $\Omega_{\nu}^{1}$ to  $\Omega_{\nu}^{\gamma}$
and applying exact closures to formulate the equations for the unknowns on  $\Omega_{\nu}^{\gamma}$.
The LOQD equations on the coarse grid $\Omega_{\nu}^{\gamma}$  are given by
\begin{subequations}\label{gloqd-gamma}
\begin{equation}\label{gloqd-1-gamma}
\frac{\partial E_p^{\gamma}   }{\partial t}
+ \frac{\partial F_p^{\gamma}   }{\partial x}
+ c \bar \sigma_{E,p}^{\gamma} E_p^{\gamma}
=  2  \bar \sigma_{B, p}^{\gamma} B_p^{\gamma} \, ,
\end{equation}
\begin{equation}\label{gloqd-2-gamma}
\frac{1}{c}\frac{\partial F_p^{\gamma}   }{\partial t}
+ c\frac{\partial (\bar f_p^{\gamma} E_p^{\gamma})   }{\partial x}
+ \bar \sigma_{R,p}^{\gamma}  F_p^{\gamma} + \bar \eta_p^{\gamma} E_p^{\gamma} = 0 \,
\end{equation}
\end{subequations}
with the boundary conditions
\begin{subequations}\label{gloqd-bcs-gamma}
\begin{equation}
F_p^{\gamma} \big|_{x=0}  = \big( c \, \bar C_p ^{\gamma-} (E_p^{\gamma}  - E_p ^{\gamma,in+}) +   F_p^{\gamma,in+} \big) \big|_{x=0}\, ,
\end{equation}
\begin{equation}
F_p^{\gamma} \big|_{x=X} = \big( c  \, \bar C_p ^{\gamma+}  (E_p^{\gamma} - E_p^{\gamma,in-}) +   F_p^{\gamma,in-} \big) \big|_{x=X} \,
\end{equation}
\end{subequations}
and the initial conditions
\begin{equation}\label{gloqd-ic-gamma}
E_p^{\gamma}\big|_{t=0 }  = E_p^{\gamma,0} \, , \quad  F_p^{\gamma} \big|_{t=0 }  = F_p^{\gamma,0} \, .
 \end{equation}
 Hereafter the equations \eqref{gloqd-gamma} are referred to as coarse-group equations.
 The coefficients of Eqs. \eqref{gloqd-gamma} are averaged with the solution on  $\Omega_{\nu}^{1}$ and defined as follows:
\begin{equation}  \label{coeff-gamma}
\bar f_p^{\gamma} = \big<  f \big>_{E,p}^{\gamma} \, , \
\bar C_p^{\gamma \pm}  = \big<  C^{\pm} \big>_{E,p}^{\gamma}  \, ,   \
\bar \sigma_{E,p}^{\gamma} = \big<  \sigma_E \big>_{E,p}^{\gamma} \, , \
\bar \sigma_{B,p}^{\gamma} = \big<  \sigma_B \big>_{B,p}^{\gamma} \, , \
\bar \sigma_{R,p}^{\gamma} = \big< \sigma_R \big>_{|F|,p}^{\gamma} \, ,
\end{equation}
where the notations for the averaged functions are given by
\begin{equation} \label{g-notation}
\big<  \psi \big>_{H,p}^{\gamma} =  \frac{\displaystyle{\sum_{g \in \tilde \Lambda_p^{\gamma} }}    \psi_{g}^1 H_{g}^1}
{\displaystyle{\sum_{g \in \tilde \Lambda_p^{\gamma} }} H_{g}^1}   \, ,
\quad
H_{g}^1 = \begin{cases}
E_{g}^1  \ \mbox{for} \   H=E  \, , \\
B_{g}^1 \ \mbox{for} \ H=B  \, ,  \\
|F_{g}^1| \ \mbox{for}  \ H=|F|    \, .
\end{cases}
\end{equation}
The compensation term in the first moment equation \eqref{gloqd-2-gamma}  is defined by
  \begin{equation} \label{eta-gamma}
\bar \eta_p^{\gamma}  =
 \frac{\displaystyle{\sum_{g \in \tilde \Lambda_p^{\gamma} }}  ( \sigma_{R,g}^1 - \bar \sigma_{R,p}^{\gamma}) F_{g}^1}
{\displaystyle{\sum_{g \in \tilde \Lambda_p^{\gamma} }}  E_{g}^1} \, .
\end{equation}
The operator form of the  LOQD equations on the  grid $\Omega_{\nu}^{\gamma}$  is  given by
\begin{equation} \label{coarse-loqd-gamma}
 \frac{\partial \mathbf{Y}_p^{\gamma}}{\partial t} + \mathcal{M}_p^{\gamma} \mathbf{Y}_p^{\gamma}   =  \mathbf{Q}_p^{\gamma} \, ,
 \
\mathbf{Y}_p^{\gamma}  =  ( E_p^{\gamma},  F_p^{\gamma})^T \, ,
\quad
p \in \mathbb{N}(n_{\nu}^{\gamma}) \, .
\end{equation}

The solution  on the  grid $\Omega_{\nu}^{\gamma}$  is used to form the grey LOQD equations on
the grid $\Omega_{\nu}^{\Gamma}$  that are defined by
\begin{equation} \label{coarse-loqd-Gamma}
 \frac{\partial \mathbf{Y}_1^{\Gamma}}{\partial t} + \mathcal{M}_1^{\Gamma} \mathbf{Y}_1^{\Gamma}   =  \mathbf{Q}_1^{\Gamma} \, ,
 \
\mathbf{Y}_1^{\Gamma}  =  ( E_1^{\Gamma},  F_1^{\Gamma})^T \, .
\end{equation}
The coefficients of Eq. \eqref{coarse-loqd-Gamma} are given by
\begin{equation} \label{coeff-Gamma}
\bar f_1^{\Gamma} = \big<  f \big>_{E}^{\gamma \to \Gamma} \, , \
\bar C_1^{\Gamma \pm}  = \big<  C^{\pm} \big>_{E}^{\gamma \to \Gamma}  \, ,   \
\bar \sigma_{E,1}^{\Gamma} = \big<  \sigma_E \big>_{E}^{\gamma \to \Gamma} \, , \
\bar \sigma_{B,1}^{\Gamma} = \big<  \sigma_B \big>_{B}^{\gamma \to \Gamma} \, , \
\bar \sigma_{R,1}^{\Gamma} = \big< \sigma_R \big>_{|F|}^{\gamma \to \Gamma} \, ,
\end{equation}
  \begin{equation}\label{eta-Gamma}
\bar \eta_1^{\Gamma}  =
 \frac{\displaystyle{\sum_{p=1}^{n_{\nu}^{\gamma}} }  (\bar \sigma_{R,p}^{\gamma} - \bar \sigma_{R,1}^{\Gamma}) F_{p}^{\gamma}}
{\displaystyle{\sum_{p=1}^{n_{\nu}^{\gamma}}  E_{p}^{\gamma}} } \, ,
\end{equation}
where
\begin{equation} \label{g2G-notation}
\big< \psi \big>_{H}^{\gamma \to \Gamma} =
 \frac{\displaystyle{\sum_{p=1}^{n_{\nu}^{\gamma}} }   \bar \psi_{p}^{\gamma} H_{p}^{\gamma}}
{\displaystyle{\sum_{p=1}^{n_{\nu}^{\gamma}} H_{p}^{\gamma}} }  \, ,
\quad
H_{p}^{\gamma} = \begin{cases}
E_{p}^{\gamma}  \ \mbox{for} \   H=E  \, , \\
B_{p}^{\gamma} \ \mbox{for} \ H=B  \, ,  \\
|F_{p}^{\gamma}| \ \mbox{for}  \ H=|F|    \, .
\end{cases}
\end{equation}
The MEB equation   coupled with the LOQD equations on the coarsest grid
$\Omega^{\Gamma}_{\nu}$     has the following form:
\begin{equation}\label{eb-Gamma}
\frac{\partial \varepsilon(T)}{\partial t} =  c  \big( \bar \sigma_{E,1}^{\Gamma} E_1^{\Gamma}  - \bar\sigma_{B,1}^{\Gamma} a_R T^4\big) \, .\
\end{equation}
The MLQD method  for TRT problems on  the sequence of   grids in  photon frequency $\mathcal{G}_{\Omega}^{\Gamma}$
 is defined by Eqs. \eqref{rt},  \eqref{mloqd}, \eqref{gloqd-gamma},   \eqref{coarse-loqd-Gamma}, and \eqref{eb-Gamma}.

\begin{figure}[th!]
  \centering
   \subfloat[ $W$-cycle, $\Gamma=3$ \label{W-cycle}]{
   \setlength{\unitlength}{0.5cm}
   \begin{picture}(6.5,5)
\put(0,5){\makebox(0,0)[c]{\tiny $\gamma$} }
\put(0,4){\makebox(0,0)[c]{\tiny \textsf{1}} }
\put(0,2){\makebox(0,0)[c]{\tiny \textsf{2}} }
\put(0,0){\makebox(0,0)[c]{\tiny \textsf{3}} }
\put(2,4){\line(1.5,-4){1.5}}
\put(3.5,0){\line(1.5,4){0.75}}
\put(4.25,2){\line(1.5,-4){0.75}}
\put(5,0){\line(1.5,4){1.5}}
\put(2,4){\circle*{0.3}}
\put(1.5,4){\makebox(0,0)[c]{\tiny \textsf{E}$_g^1$} }
\put(3.5,0){\circle*{0.3}}
\put(3.1,0){\makebox(0,0)[c]{\tiny \textsf{T}} }
\put(4.25,2.){\circle*{0.3}}
\put(3.75,2){\makebox(0,0)[c]{\tiny \textsf{E}$_p^{\gamma}$} }
\put(5,0){\circle*{0.3}}
\put(4.6,0){\makebox(0,0)[c]{\tiny \textsf{T}} }
\end{picture}
}
\hspace{1cm}
   \subfloat[ $F$-cycle, $\Gamma=4$ \label{F-cycle}]{
   \setlength{\unitlength}{0.5cm}
   \begin{picture}(8,5)
\put(0,5){\makebox(0,0)[c]{\tiny $\gamma$} }
\put(0,4){\makebox(0,0)[c]{\tiny \textsf{1}} }
\put(0,2.66){\makebox(0,0)[c]{\tiny \textsf{2}} }
\put(0,1.33){\makebox(0,0)[c]{\tiny \textsf{3}} }
\put(0,0){\makebox(0,0)[c]{\tiny \textsf{4}} }
\put(2,4){\line(1.5,-4){1.5}}
\put(3.5,0){\line(1.5,4){0.5}}
\put(4,1.33){\line(1.5,-4){0.5}}
\put(4.5,0){\line(1.5,4){1.}}
\put(5.5,2.66){\line(1.5,-4){1.}}
\put(6.5,0){\line(1.5,4){1.5}}
\put(2,4){\circle*{0.3}}
\put(1.5,4){\makebox(0,0)[c]{\tiny \textsf{E}$_g^1$} }
\put(3.5,0){\circle*{0.3}}
\put(3.1,0){\makebox(0,0)[c]{\tiny \textsf{T}} }
\put(4.,1.33){\circle*{0.3}}
\put(3.5,1.33){\makebox(0,0)[c]{\tiny \textsf{E}$_p^{\gamma}$} }
\put(4.5,0){\circle*{0.3}}
\put(4.1,0){\makebox(0,0)[c]{\tiny \textsf{T}} }
\put(5.5,2.66){\circle*{0.3}}
\put(5,2.66){\makebox(0,0)[c]{\tiny \textsf{E}$_p^{\gamma}$} }
\put(6.5,0){\circle*{0.3}}
\put(6.1,0){\makebox(0,0)[c]{\tiny \textsf{T}} }
\end{picture}}
\hspace{1cm}
\subfloat[ $V$-cycle, $\Gamma=2$ \label{V-cycle}]{
\setlength{\unitlength}{0.5cm}
\begin{picture}(5,5)
\put(0,5){\makebox(0,0)[c]{\tiny $\gamma$} }
\put(0,4){\makebox(0,0)[c]{\tiny \textsf{1}} }
\put(0,0){\makebox(0,0)[c]{\tiny \textsf{2}} }
\put(2,4){\line(1.5,-4){1.5}}
\put(3.5,0){\line(1.5,4){1.5}}
\put(2,4){\circle*{0.3}}
\put(1.5,4){\makebox(0,0)[c]{\tiny \textsf{E}$_g^1$} }
\put(3.5,0){\circle*{0.3}}
\put(3.1,0){\makebox(0,0)[c]{\tiny \textsf{T}} }
\end{picture}
}
   \caption{Diagrams of
  multigrid      cycles on  hierarchies of $\Gamma$ grids in frequency.
  \textsf{E$_g^1$} - calculation  of spectrum on the fine  grid  $\Omega_{\nu}^{1}$  by solving the multigroup LOQD equations,
 \textsf{T} - calculation of temperature by solving coupled  grey LOQD and MEB equations on $\Omega_{\nu}^{\Gamma}$,
  \textsf{E$_p^{\gamma}$} - calculation of spectrum on $\Omega_{\nu}^{\gamma}$ ($\gamma>1$) by solving  the coarse-group LOQD equations.}
  \label{cycles}
\end{figure}
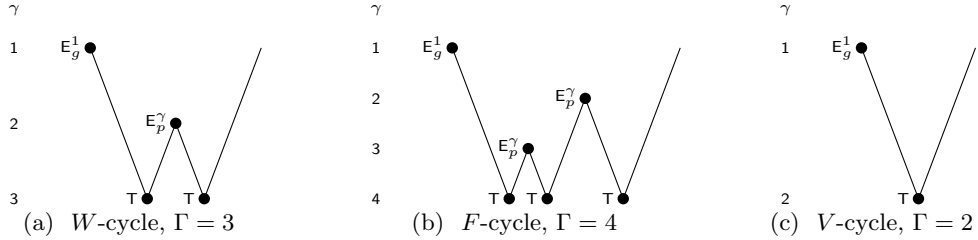

To solve the hierarchy of  LOQD equations coupled with the  MEB equation,
we apply   $W$  and full ($F$)  multigrid cycles  to visit grids   \cite{multigrid}.
 Algorithm \ref{mlqd-multigrid}  describes the MLQD method on the hierarchy of frequency grids $\mathcal{G}_{\Omega}^{\Gamma}$.
 The set  $\mathcal{S}=\{\gamma_k, k \in \mathbb{N}(K)\}$  defines
   the schedule of visiting grids  after each evaluation of temperature on the coarsest grid $\Omega^{\Gamma}_{\nu}$.
 Figure \ref{W-cycle}   shows  the diagram  of  the $W$-cycle for three grids ($\Gamma=3$)
for which  $\mathcal{S} = \{2\}$ and $K=1$.
The $F$-cycle   on four  grids ($\Gamma=4$)
 with $K=2$ and the schedule  $\mathcal{S} = \{3,2\}$
 is illustrated in   Fig.   \ref{F-cycle}.
 The MLQD method with two grids (Sec. \ref{sec:mlqd-2grid}) is equivalent  to Algorithm \ref{mlqd-multigrid}  with the $V$-cycle  the diagram of which is shown in Fig. \ref{V-cycle}.
 The cycles are executed until either the convenience criteria  for $T$ and $E$ are satisfied or the number of cycles reaches
 the given maximum number $\ell_{max}$.

The multigrid cycle  starts  from calculating the spectrum on the fine grid $\Omega_{\nu}^1$  for
 the current estimation of temperature. The new temperature
 is evaluated by solving the effective grey problem on $\Omega^{\Gamma}_{\nu}$  formulated with the fine-grid spectrum.
 Then the spectrum is computed on some  coarse grid  $\Omega_{\nu}^{\gamma}$ for $1<\gamma < \Gamma$
 using the new temperature. On the next stage, the obtained coarse-grid spectrum on $\Omega_{\nu}^{\gamma}$  is used to
 average opacities and factors and   calculate their grey quantities. This forms the updated effective grey problem on $\Omega^{\Gamma}_{\nu}$
 that is solved to get the   temperature  at this  stage of the cycle.
 Then these two elements of the cycle are repeated using the coarse grid according to the schedule  $\mathcal{S}$.

The coarse-group LOQD equations on $\Omega_{\nu}^{\gamma}$  are derived by exact averaging
the LOQD equations on $\Omega_{\nu}^{\gamma-1}$. As a results,
if the group LOQD equations on these  grids are defined with the same $T$, then
\begin{equation}
E_p^{\gamma} = \sum_{p' \in \Lambda_p^{\gamma} } E_{p'}^{\gamma-1} =
\sum_{g \in \tilde \Lambda_p^{\gamma} } E_g^1 \, ,
\quad \quad
F_p^{\gamma} =
\sum_{p' \in \Lambda_p^{\gamma} } F_{p'}^{\gamma-1} =
\sum_{g \in \tilde \Lambda_p^{\gamma} } F_g^1\, .
\end{equation}
 Thus, there is no change in spectrum from grid to grid  in this case.
The multigrid cycles are defined in such a way that
 a grid is visited only after the  temperature  update was performed on the coarsest grid $\Omega_{\nu}^{\Gamma}$.
If the  spectrum is evaluated  on $\Omega_{\nu}^{\gamma'}$  with  the current estimate of temperature
then all coarser grids with  $\gamma' < \gamma < \Gamma$ are skipped on the way down to $\Omega_{\nu}^{\Gamma}$.
After evaluation of temperature on $\Omega_{\nu}^{\Gamma}$  the algorithm is scheduled to update the
spectrum on the grid $\Omega_{\nu}^{\gamma''}$.  The algorithm skips all grids with $ \gamma''  < \gamma  < \Gamma$
on the way up to the scheduled finer grid.

\vspace{1cm}
\begin{algorithm}[H]
\DontPrintSemicolon
$s=-1$, $T^{(0)}=T^{j-1}$, $f_{g}^{(0)}=f_{g}^{j-1}$\;
\While{$|| T^{(s)} - T^{(s-1)} ||>\epsilon||T^{(s)}|| \, \& \,
 ||E^{(s)} - E^{(s-1)}||> \epsilon||E^{(s)}|| $ }{
$\bullet$ transport (outer) iterations\;
 $s=s+1$ \;
  \If{$s>0$}{
$T^{(s)} \Rightarrow \mathcal{L}_g^{(s)}$\;
$c^{-1}  \partial_t I_g  + \mathcal{L}_g^{(s)} I_g = q_g[T^{(s)}]   \ \mbox{on} \  \Omega_{\nu}^1 \Rightarrow I_g^{(s)}, g \in \mathbb{N}(n_{\nu}^1)$\;
 $I_g^{(s)} \Rightarrow f_{g}^{(s)}, C_g^{\pm (s)}, g \in \mathbb{N}(n_{\nu}^1)$\;
 }
 $\ell=-1$,  $T^{(1,s)}=T^{(s)}$\;
\While{$\ell \le \ell_{max}$ or  $|| T^{(\ell,s)} - T^{(\ell-1,s)} ||>\tilde\epsilon||T^{(\ell,s)}||  \, \& \,
 || E^{(\ell,s)} - E^{(\ell-1,s)} ||>\tilde\epsilon||E^{(\ell,s)}||$  }{
 $\bullet$ low-order (inner) iterations\;
 $\ell=\ell+1$ \;
  calculation of fine-grid    spectrum on $\Omega_{\nu}^{1}$\;
$T^{(\ell,s)}, f_g^{(s)},C_g^{\pm (s)} \Rightarrow \mathcal{M}_g^{1(\ell,s)}$\;
$\partial _t \mathbf{Y}_g^1+ \mathcal{M}_g ^{1(\ell,s)}\mathbf{Y}_g^1 = \mathbf{Q}_g^1[T^{(\ell,s)}]    \  \mbox{on} \ \Omega_{\nu}^1 \Rightarrow E_g^{1(\ell,s)}, F_g^{1(\ell,s)}, g \in \mathbb{N}(n_{\nu}^1)$ \;
 $\tilde T^{[1]} = T^{(\ell,s)}$ \;
 \For{$k \leftarrow 0  \ \KwTo  \ K $}{
 $\gamma \leftarrow 1$\;
 \If{$k>0$}{
 $\gamma  \leftarrow \gamma_k \in \mathcal{S} $\;
calculation of coarse-grid   spectrum on $\Omega_{\nu}^{\gamma}$\;
$\tilde T^{[k]}, E_g^{1(\ell,s)}, F_g^{1(\ell,s)}, f_g^{(s)}, C_g^{\pm (s)} \Rightarrow \mathcal{M}_p^{\gamma[k]}$\;
$\partial _t \mathbf{Y}_p^{\gamma}+ \mathcal{M}_p ^{\gamma[k]}\mathbf{Y}_p^{\gamma} = \mathbf{Q}_p[\tilde T^{[k]}]    \  \mbox{on} \ \Omega_{\nu}^{\gamma}\Rightarrow E_p^{\gamma[k]}, F_p^{\gamma[k]}, p \in \mathbb{N}(n_{\nu}^{\gamma})$ \;
 }
 $\tilde T^{[k]}, E_p^{\gamma[k]}, F_p^{\gamma[k]} \Rightarrow \mathcal{M}_1^{\Gamma[k]}$\;
$\partial _t \mathbf{Y}_1^{\Gamma} + \mathcal{M}_1^{\Gamma[k]} \mathbf{Y}_1^{\Gamma} = \mathbf{Q}_1^{\Gamma}$ on $\Omega_{\nu}^{\Gamma}$  \& MEB Eq. \eqref{eb-Gamma} $\Rightarrow \tilde T^{[k+1]}, E_1^{\Gamma}, F_1^{\Gamma}$\;
}
$ T^{(\ell+1,s)} \leftarrow \tilde T^{[k+1]}$\;
}
$T^{(s+1)} \leftarrow  T^{(\ell+1,s)}$ \;
}
\caption{\label{mlqd-multigrid}  The MLQD method with multiple  grids in frequency
 $\Omega_{\nu}^{\gamma}$, $\gamma \in \mathbb{N}(\Gamma)$.}
\end{algorithm}

\section{\label{sec:disc} Discretization of Equations}

The system of the high-order RT,  LOQD, and MEB equations
is approximated by the  implicit Euler time integration scheme.
The opacities and emission terms are evaluated at the current time level.
Thus, the implicitly balanced  time integration scheme is applied \cite{dana-2007}.
The RT equation is approximated in space with the simple corner balance method \cite{mla-ttsp-1997}.
The second-order finite volume  method is applied
to discretize  the  multigroup LOQD equations on  $\Omega_{\nu}^1 $ over space \cite{dya-jcp-2019}.
We define the spatial mesh
$\{ x_{i-1/2},    i~\in~\mathbb{N}(n_x~+~1), x_{1/2}~=~0,  x_{n_x+1/2}=X \}$.
The photon balance equation \eqref{mloqd-1} is integrated over the $i$-th cell
 ($x_{i-1/2}\le x \le x_{i+1/2}$). The  first moment equation  \eqref{mloqd-2}
 is integrated over  $x_{i-1} \le x \le x_{i}$, where $x_i=\frac{1}{2} (x_{i-1/2} + x_{i+1/2})$.
 The discretized LOQD equations  \eqref{mloqd} on $\Omega_{\nu}^1$  at $t=t^j$ have the following form:
 \begin{subequations}\label{mloqd-d}
\begin{equation} \label{mloqd1-d}
\frac{\Delta x_i}{\Delta t^{j}} \Big( E_{g,i}^{1,j} - E_{g,i}^{1,j-1} \Big)
+F_{g,i+1/2}^{1,j} - F_{g,i-1/2}^{1,j}
+ c \sigma_{E,g,i}^{1,j} \Delta x_j E_{g,i}^{1,j} =
\sigma_{B,g,i}^{1,j}\Delta x_i B_{g,i}^{1,j}  \, ,
\end{equation}
\begin{multline} \label{mloqd2-d}
   \frac{\Delta x_{i+1/2}}{c\Delta t^{j}}    \Big(F_{g,i + 1/2}^{1,j} -   F_{g,i + 1/2}^{1,j-1} \Big)
+  c \Big(f_{g,i+1}^{1,j} E_{g,i+1}^{1,j} -   f_{g,i}^{1,j} E_{g,i}^{1,j}\Big) \\
+  \sigma_{R,g,i+1/2}^{1,j} \Delta x_{i+1/2} F_{g,i+1/2}^{1,j} = 0 \, ,
\end{multline}
 \end{subequations}
where
\begin{equation}
 \Delta x_i = x_{i+1/2} - x_{i-1/2} \, ,
\end{equation}
\begin{equation}
\sigma_{R,g,i+1/2}^{1,j} = \frac{\sigma_{R,g,i}^{1,j} \Delta x_i   +  \sigma_{R,g,i+1}^{1,j} \Delta x_{i+1}}{\Delta x_i+ \Delta x_{i+1}} \, .
\end{equation}
 $j$ is the index of the time step.
Integer $\pm \frac{1}{2}$ subscripts refer to cell-edge quantities, and integer subscripts refer to cell-average
quantities.

The spatial approximation of the coarse-group LOQD equations on $\Omega_{\nu}^{\gamma}$ is algebraically consistent with the group LOQD equations on the fine grid $\Omega_{\nu}^{1}$. The discretized  LOQD equations on $\Omega_{\nu}^{\gamma}$  are given by
 \begin{subequations}\label{mloqd-d-gamma}
\begin{equation} \label{mloqd1-d-gamma}
\frac{\Delta x_i}{\Delta t^{j}} \Big( E_{p,i}^{\gamma,j} - E_{p,i}^{\gamma,j-1} \Big)
+F_{p,i+1/2}^{\gamma,j} - F_{p,i-1/2}^{\gamma,j}
+ c \bar \sigma_{E,p,i}^{\gamma,j} \Delta x_j E_{p,i}^{\gamma,j} =
\bar \sigma_{B,p,i}^{\gamma,j}\Delta x_i B_{p,i}^{\gamma,j}  \, ,
\end{equation}
\begin{multline} \label{mloqd2-d-gamma}
   \frac{\Delta x_{i+1/2}}{c\Delta t^{j}}    \Big(F_{p,i + 1/2}^{\gamma,j} -   F_{p,i + 1/2}^{\gamma,j-1} \Big)
+  c \Big(   \big(\bar f_{p,i+1}^{\gamma,j} +  \hat \eta_{p,i+1/2}^{\gamma,j} \big) E_{p,i+1}^{\gamma,j}
 -   \big(\bar f_{p,i}^{j} +  \check \eta_{p,i+1/2}^{\gamma,j} \big)E_{p,i}^{\gamma,j}\Big) \\
+  \bar \sigma_{R,p,i+1/2}^{\gamma,j} \Delta x_{i+1/2} F_{p,i+1/2}^{\gamma,j} = 0 \, ,
\end{multline}
 \end{subequations}
where
\begin{equation}
\bar f_{p,i}^{\gamma,j} = \big<  f_i^j \big>_{E,p}^{\gamma} \, , \
\bar \sigma_{E,p,i}^{\gamma,j}  = \big<  \sigma_{E,i}^j \big>_{E,p}^{\gamma} \, , \
\bar \sigma_{B,p,i}^{\gamma,j}  = \big<  \sigma_{B,i}^j \big>_{B,p}^{\gamma} \, , \
\bar \sigma_{R,p,i+1/2}^{\gamma,j} = \big< \sigma_{R,i+1/2}^j \big>_{|F|,p}^{\gamma} \, ,
\end{equation}
\begin{equation}
\big<  \psi_{\alpha}^j \big>_{H,p}^{\gamma} =  \frac{\displaystyle{\sum_{g \in  \tilde \Lambda_p^{\gamma} }}    \psi_{g,\alpha}^{1,j} H_{g,\alpha}^{1,j}}
{\displaystyle{\sum_{g \in \tilde \Lambda_p^{\gamma} }} H_{g,\alpha}^{1,j}}   \, ,
\quad
H_{g,\alpha}^{1,j} = \begin{cases}
E_{g,i}^{1,j} \, , \alpha=i \,   \ \mbox{for} \   H=E  \, , \\
B_{g,i}^{1,j} \, , \alpha=i \,  \ \mbox{for} \ H=B  \, ,  \\
|F_{g,i+1/2}^{1,j}|  \, , \alpha=i+\frac{1}{2} \,   \ \mbox{for}  \ H=|F|    \, .
\end{cases}
\end{equation}
\begin{equation}
 \hat \eta_{p,i+1/2}^{\gamma,j} =\begin{cases}
 \frac{\xi_{p,i+1/2}^{\gamma,j}}{\displaystyle{c \sum_{g \in \tilde \Lambda_p^{\gamma} }  E_{g,i+1}^{1,j}}} \
 & \mbox{for} \  \xi_{p,i+1/2}^{\gamma,j} > 0 \, ,\\
 0   & \mbox{for} \  \xi_{p,i+1/2}^{\gamma,j} \le 0 \, ,
 \end{cases}
 \quad
  \check \eta_{p,i+1/2}^{\gamma,j} =\begin{cases}
 0   & \mbox{for} \  \xi_{p,i+1/2}^{\gamma,j} \ge 0  \, ,   \\
 -\frac{\xi_{p,i+1/2}^{\gamma,j}}{\displaystyle{c \sum_{g \in \tilde \Lambda_p^{\gamma} }  E_{g,i}^{1,j}}} \
 & \mbox{for} \  \xi_{p,i+1/2}^{\gamma,j} < 0  \, ,
 \end{cases}
\end{equation}
\begin{equation}
\xi_{p,i+1/2}^{\gamma,j}  = \sum_{g \in \tilde \Lambda_p^{\gamma} }
 \Big(\sigma_{R,g,i+1/2}^{1,j} - \bar \sigma_{R,p,i+1/2}^{\gamma,j} \Big)
  F_{g,i+1/2}^{1,j} \, .
\end{equation}
Similarly, the discretized LOQD equations on the  grid $\Omega_{\nu}^{\gamma}$   (Eqs. \eqref{mloqd-d-gamma})
are averaged over   groups to derive the discrete grey LOQD equations  (Eqs. \eqref{coarse-loqd-Gamma}) on
the grid $\Omega_{\nu}^{\Gamma}$  and define   their coefficients \eqref{coeff-Gamma}.
The discretized MEB equations has the form:
\begin{equation}\label{eb-Gamma-d}
\frac{\partial \varepsilon(T_i^j\big) }{\partial t} =  c  \Big( \bar \sigma_{E,1,i}^{\Gamma,j} E_{1,i}^{\Gamma,j}  - \bar\sigma_{B,1,i}^{\Gamma,j} a_R \big(T_i^j\big)^4\Big) \, .\
\end{equation}

The  grey LOQD equations    on $\Omega_{\nu}^{\Gamma}$   and the MEB equation (Eq. \eqref{eb-Gamma-d})  are solved by Newton's method.
The grey opacity
\begin{equation}
\bar \sigma_{E,1}^{\Gamma} (T)=
 \frac{\displaystyle  \sum_{p =1}^{n_{\nu}^{\gamma}} \bar \sigma_{E,p}^{\gamma}(T) E_{p}^{\gamma}(T)}
{\displaystyle  \sum_{p =1}^{n_{\nu}^{\gamma}}  E_{p}^{\gamma}(T)}  \, , \quad
\gamma \in \mathcal{S}
 \end{equation}
 depends locally on $T$  through  $\sigma_{E,g}(T)$ on $\Omega_{\nu}^1$ and
globally through the coarse-grid  $E_{p}^{\gamma}(T) $ as well as
$E_{g}^{1}(T)$ applied to compute the opacity   $\bar \sigma_{E,p}^{\gamma}$ on $\Omega_{\nu}^{\gamma}$.
The  Fr$\acute{\mbox{e}}$chet derivative  of  $\bar \sigma_{E,1}^{\Gamma}$
  is used in   the linearized equations to account for its  variation  due to change in temperature.
 The Fr$\acute{\mbox{e}}$chet derivative
 $\mathcal{D}  \bar \sigma_{E,1}^{\Gamma}$
  is  a linear operator such that
 \begin{equation}
 \bar \sigma_{E,1}^{\Gamma}(T + \Delta T) - \bar \sigma_{E,1}^{\Gamma}(T)  =
  \mathcal{D}   \bar \sigma_{E,1}^{\Gamma}  \Delta T + \rho(T, \Delta T) \, ,
 \end{equation}
 where
  \begin{equation}
 \frac{ \rho(T, \Delta T)}{|| \Delta T|| } \to 0 \quad  \mbox{as}  \quad || \Delta T|| \to 0 \, .
 \end{equation}
 In  discrete space, $\mathcal{D}   \bar \sigma_{E,1}^{\Gamma}$  is a matrix.
 There are different ways to estimate   $\mathcal{D}   \bar \sigma_{E,1}^{\Gamma}$.
A robust and  efficient variant is to approximate  it by the diagonal matrix given by \cite{dya-aristova-vya-mm1996,aristova-vya-avk-m&c1999}
\begin{equation} \label{f-derivative-d}
 \mathcal{D}   \bar \sigma_{E,1}^{\Gamma}
= \mbox{diag} \Big[\big(D \bar \sigma_{E,1}^{\Gamma}\big)_1^{[k]},
\ldots,\big(D \bar \sigma_{E,1}^{\Gamma}\big)_i^{[k]},\ldots,
\big(D \bar \sigma_{E,1}^{\Gamma}\big)_{n_x}^{[k]} \Big] \, ,
\end{equation}
where
\begin{equation}
\big(D \bar \sigma_{E,1}^{\Gamma}\big)_i^{[k]} =
 \frac{\bar \sigma_{E,1,i}^{\Gamma}\big(\tilde T^{[k]}\big) - \bar \sigma_{E,1,i}^{\Gamma}\big(\tilde T^{[k-1]}\big)}
 {\tilde T^{[k]}_i - \tilde T^{[k-1]}_i} \, .
 \end{equation}
The value of the discrete Fr$\acute{\mbox{e}}$chet derivative  \eqref{f-derivative-d} is fixed
during  Newton's iterations for solving of the linearized grey LOQD and MEB equations.

\section{\label{sec:res} Numerical Results}

\begin{figure}[h!]
\centering
\subfloat[Temperature   \label{T}]{\includegraphics[scale=0.3]{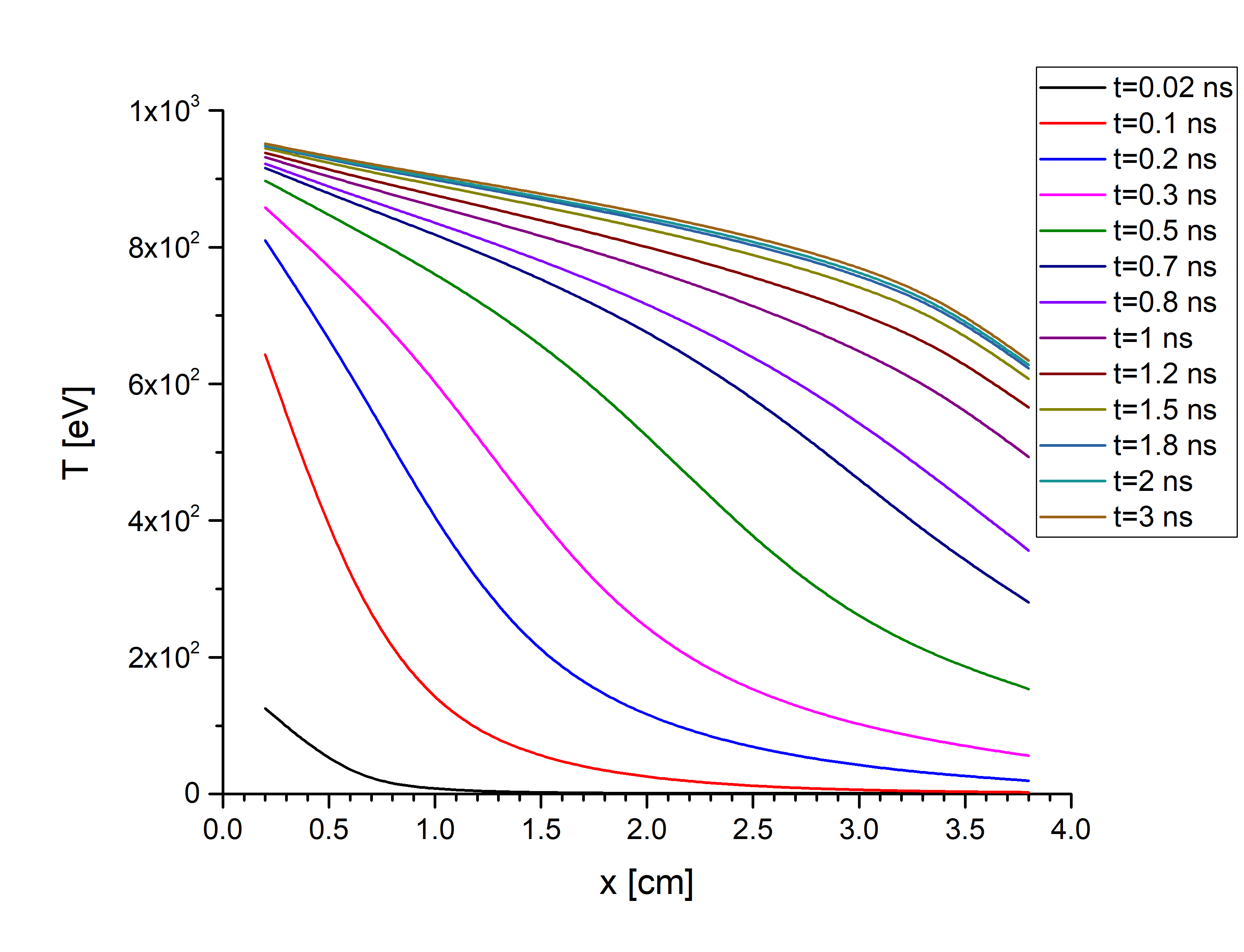}}
\subfloat[Energy density \label{E}]{\includegraphics[scale=0.3]{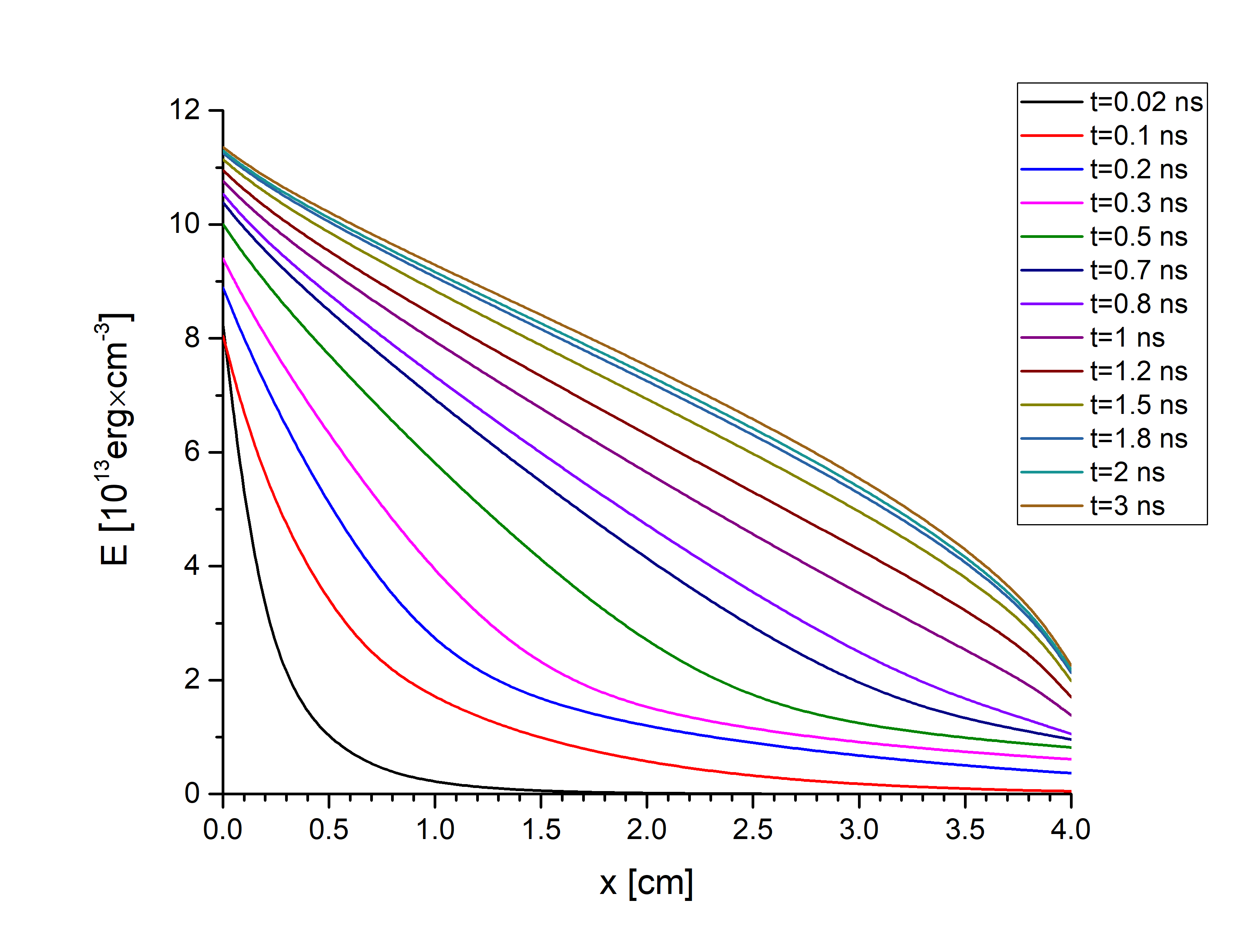}}
\caption{\label{solution}
Numerical solution of the FC test with 256 groups and $\Delta t=2 \! \times \! 10^{-2}$ ns.}
\end{figure}

In this section,  numerical  results of the Fleck-Cummings (FC) test are presented  \cite{fleck-1971}.
The test is defined for a  slab ($0\le x \le 4$ cm) with one  material.
The  spectral opacity of the material is given by
\begin{equation}
\sigma_{\nu}(T)    =    \frac{27}{(h\nu)^3}\big(1-e^{-\frac{h\nu}{kT}}\big) \, .
\end{equation}
There is incoming radiation with the black-body spectrum at $kT_b$=1 ~keV at the left boundary.
The right boundary is vacuum.
The initial temperature in the domain is $kT_0$=10$^{-3}$ keV.
\linebreak
 $ I_{g}|_{t=0}  =   B_{g}(T_0)$.
The material energy is given by $\varepsilon(T)   =  c_v T$ with  $c_v  = \    0.5917a_R  T_b^3$.
The spatial  mesh is uniform with  10 intervals.
The time interval of the problem is  $t \in [0, 3 \, \mbox{ns}]$.
The frequency grid $\Omega_{\nu}$ is defined by $n_{\nu}=256$   groups over the range $0~\le~h\nu~\le~10^7$~keV.
There are  $n_{\nu}-2$  groups evenly spaced in logarithmic scale between $h\nu_a= 10^{-4}$~keV and $h\nu_b = 10$ keV.
The double S$_8$  Gauss-Legendre quadrature set is used and hence there are 16  angular directions.
The test is calculated with  two time steps:  $\Delta t = 2 \times 10^{-2}$ ns and $\Delta t = 4 \times 10^{-2}$ ns.
The parameters of convergence criteria  are $\epsilon  = 10^{-6}$ and $\tilde \epsilon  = 10^{-7}$.
Figure \ref{solution} shows temperature and total energy density at various instants of time obtained with $\Delta t = 2 \times 10^{-2}$ ns and  illustrates evolution of    heat and radiation waves.

The algorithms with $W$ and $F$ cycles  (see Fig. \ref{cycles}) were studied  on a variety of grid  hierarchies
with different $\Gamma$.
The  grids are formed by successive  and  uniform  coarsening.
They are  defined by the number of groups, i.e. $n_{\nu}^{\gamma}$.
The number of Newton's iterations for solving the grey LOQD  and MEB equations to evaluate temperature is equal to 1.
The algorithms use different values of the parameter $\ell_{max}$ that restricts the number of inner (low-order)  iterations
and hence the number of multigrid cycles on each transport iteration.

\begin{figure}[h!]
\centering
\subfloat[$\Gamma$=2,  $n_{\nu}^{\gamma}$=256,1, $V$-cycle, $\ell_{max}$=4. \label{V-256-002-fig}]
{\includegraphics[scale=0.28]{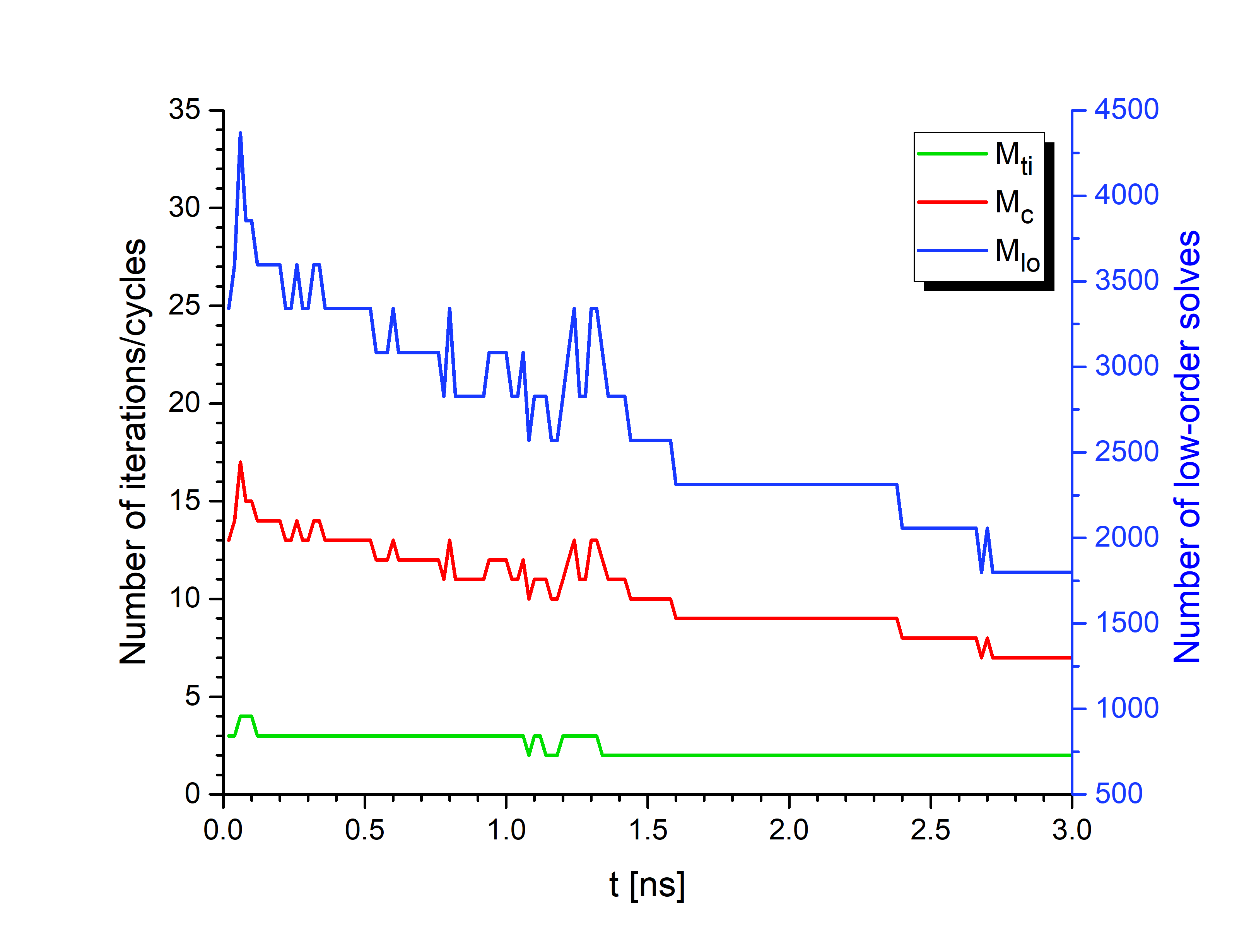}}
\hfill
\subfloat[$\Gamma$=3, $n_{\nu}^{\gamma}$=256,32,1,  $W$-cycle, $\ell_{max}$=2. \label{W-256-32-002-fig}]{\includegraphics[scale=0.28]{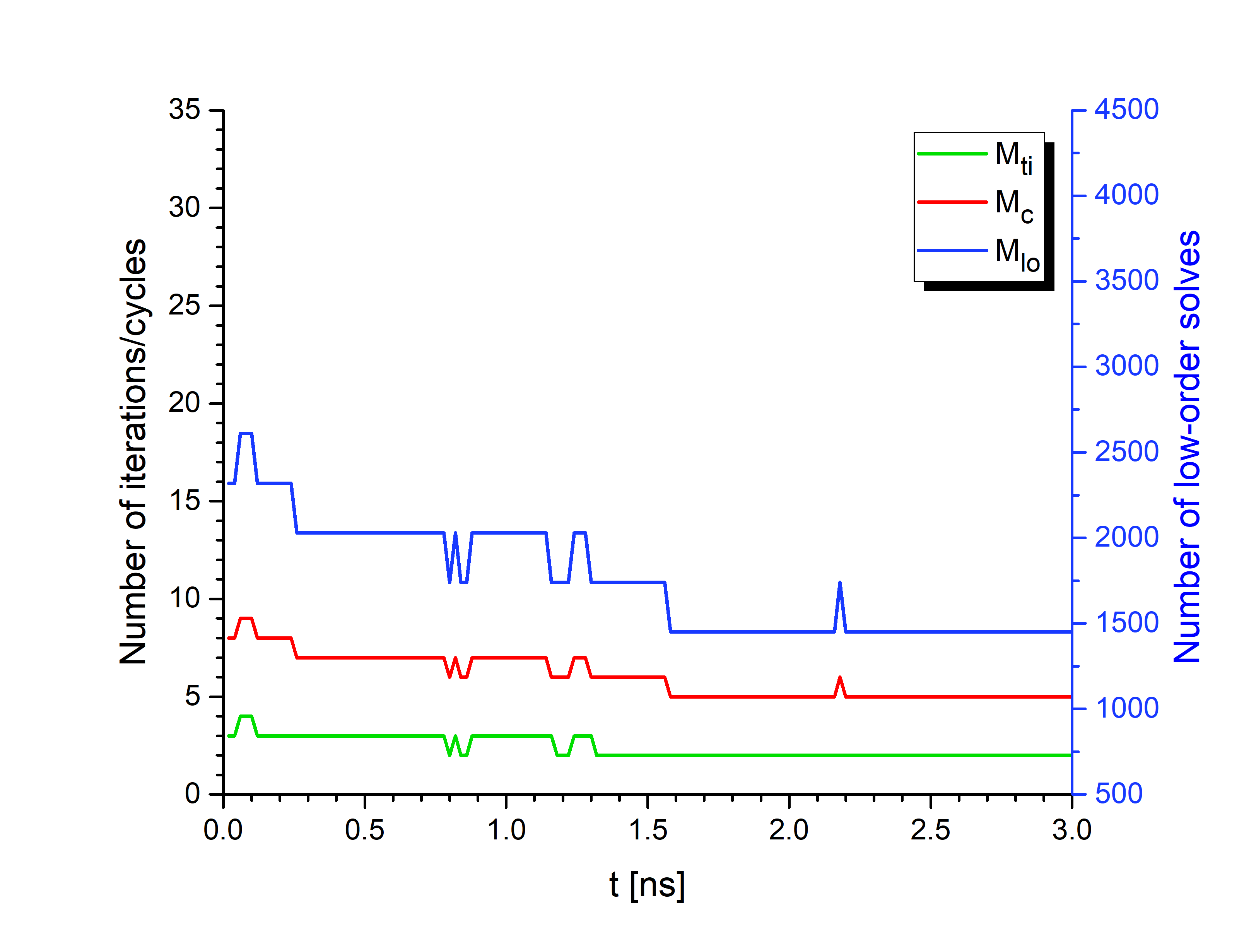}}\\
\subfloat[$\Gamma$=4,  $n_{\nu}^{\gamma}$=256,32,16,1, $F$-cycle, $\ell_{max}$=2.
\label{F-256-32-16-002-fig}]
{\includegraphics[scale=0.28]{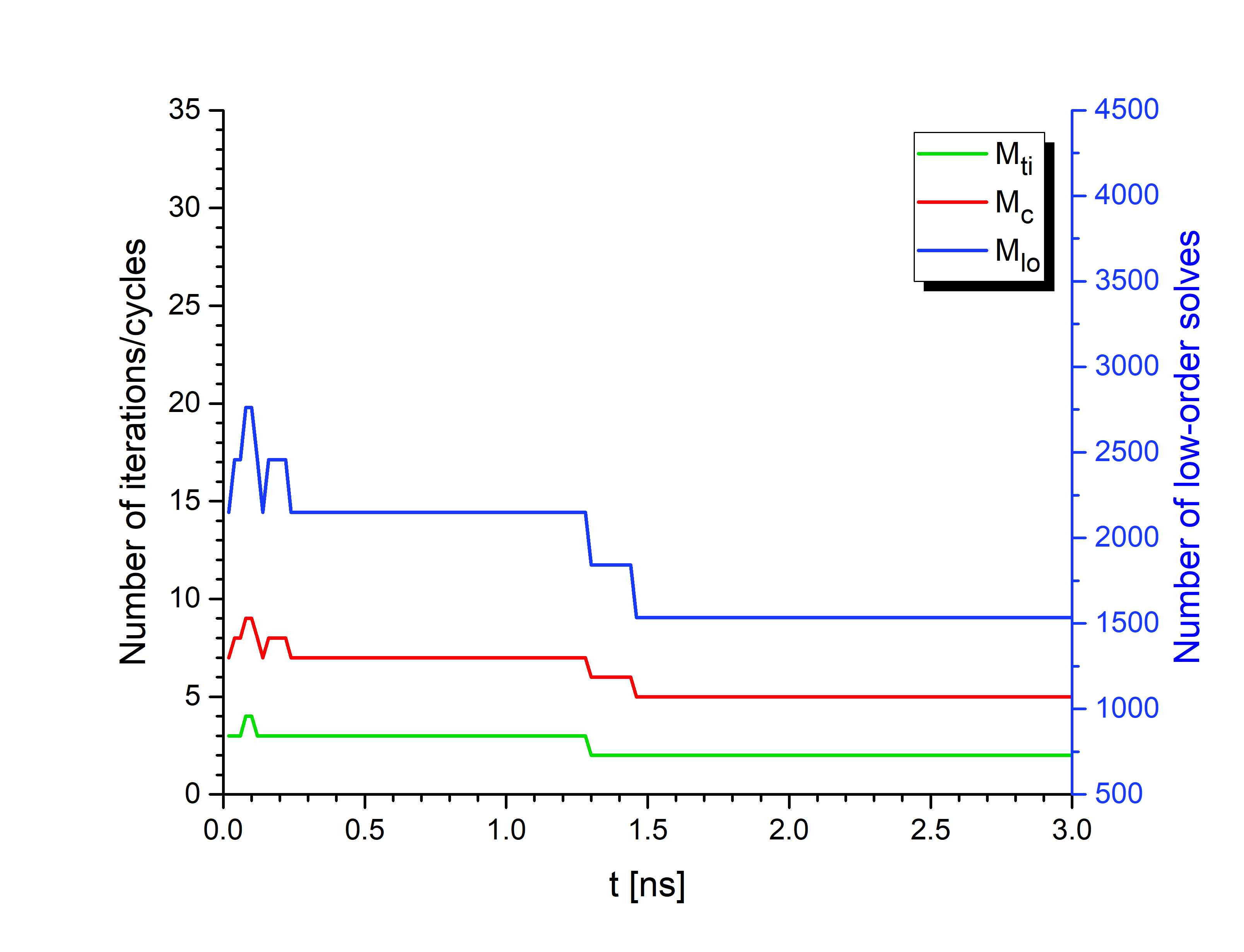}}
\hfill
\subfloat[$\Gamma$=5, $n_{\nu}^{\gamma}$=256,32,16,4,1, $F$-cycle, $\ell_{max}$=2. \label{F-256-32-16-4-002-fig}]{\includegraphics[scale=0.28]{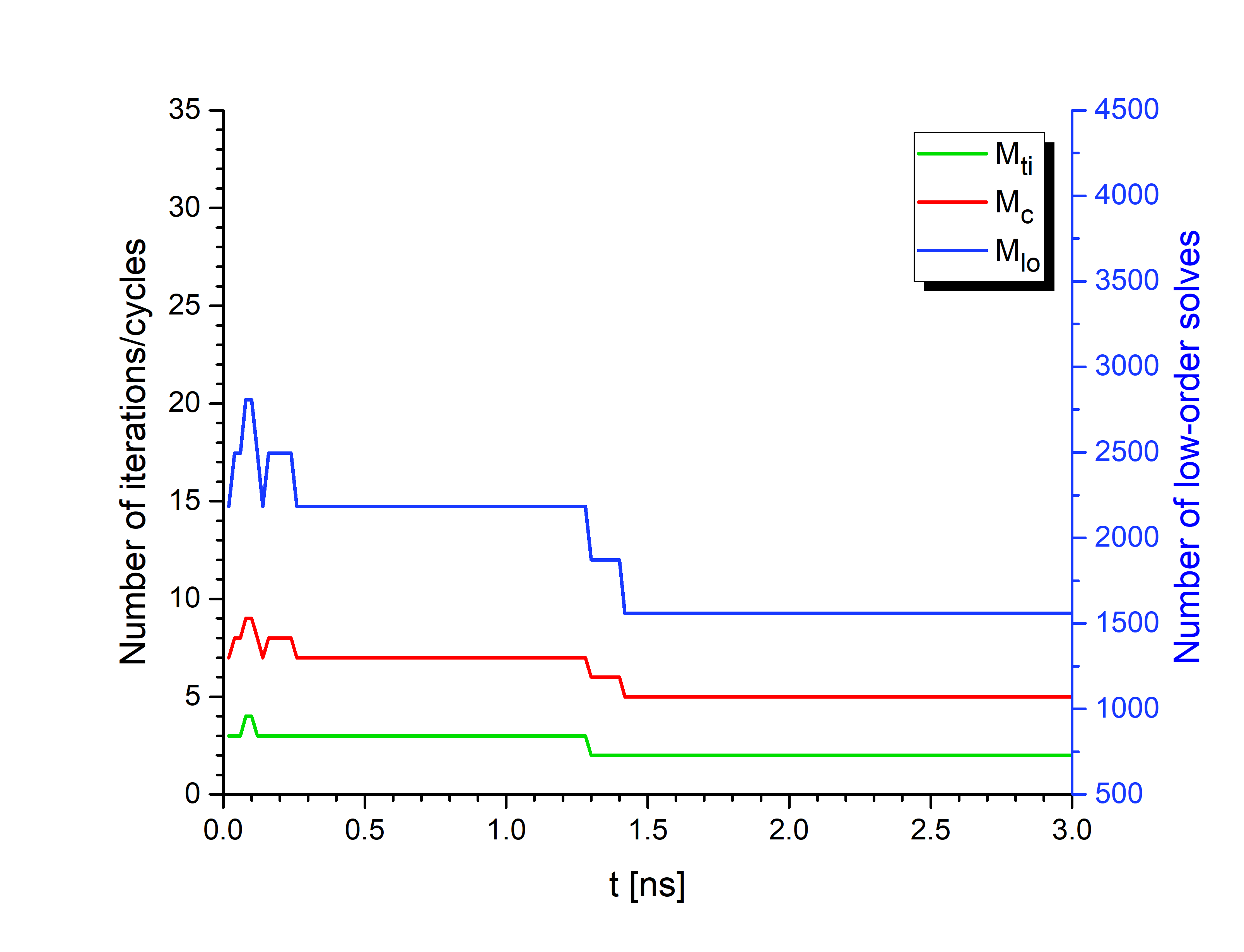}}\\
\subfloat[$\Gamma$=6, $n_{\nu}^{\gamma}$=256,128,64,32,1,  $F$-cycle, $\ell_{max}$=1. \label{F-256-128-64-32-002-fig}]{\includegraphics[scale=0.28]{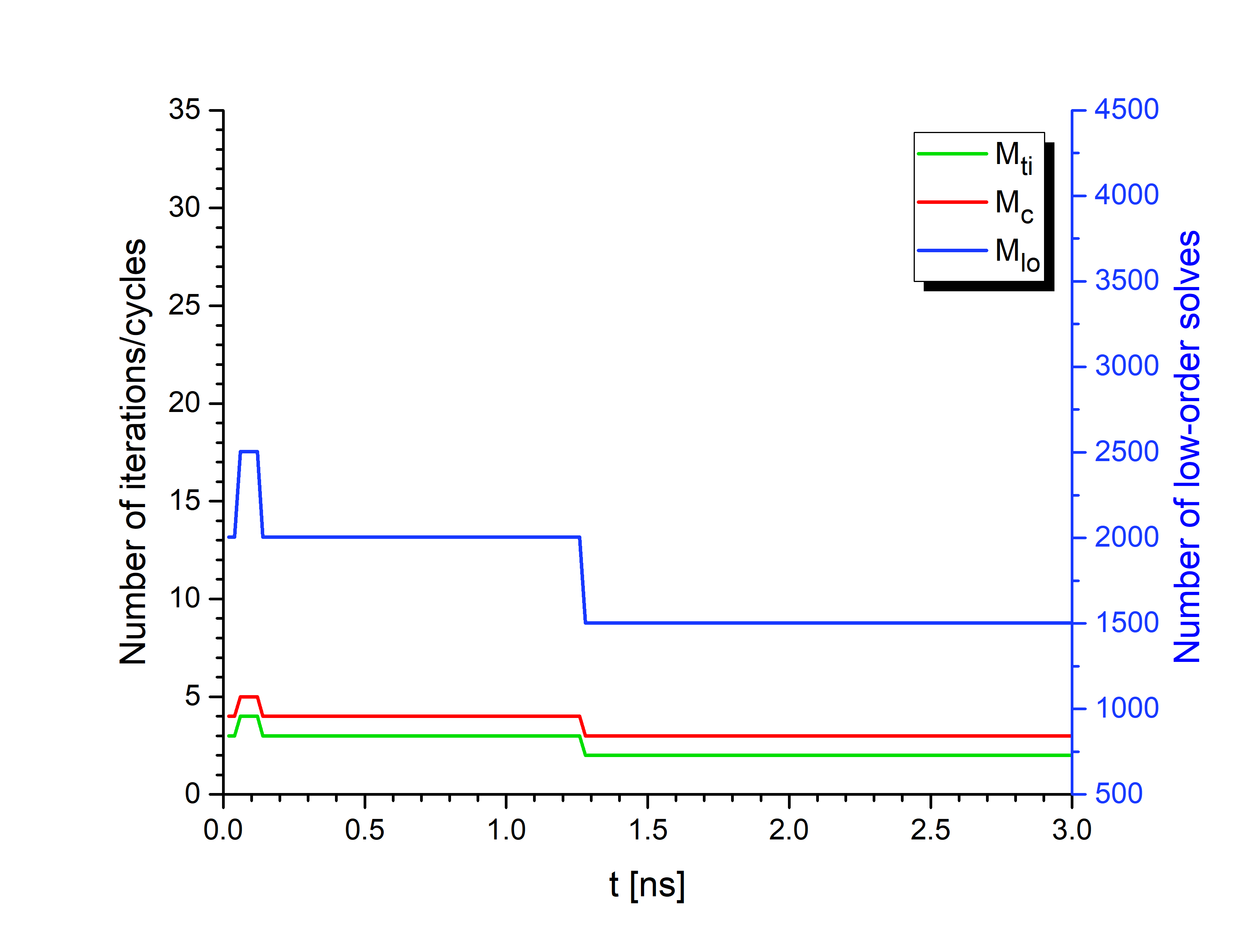}}
\hfill
\subfloat[$\Gamma$=7,\,$n_{\nu}^{\gamma}$=256,128,32,16,8,4,1,\,$F$-cycle,\,$\ell_{max}$=1. \label{F-256-128-32-16-8-4-002-fig}]{\includegraphics[scale=0.28]{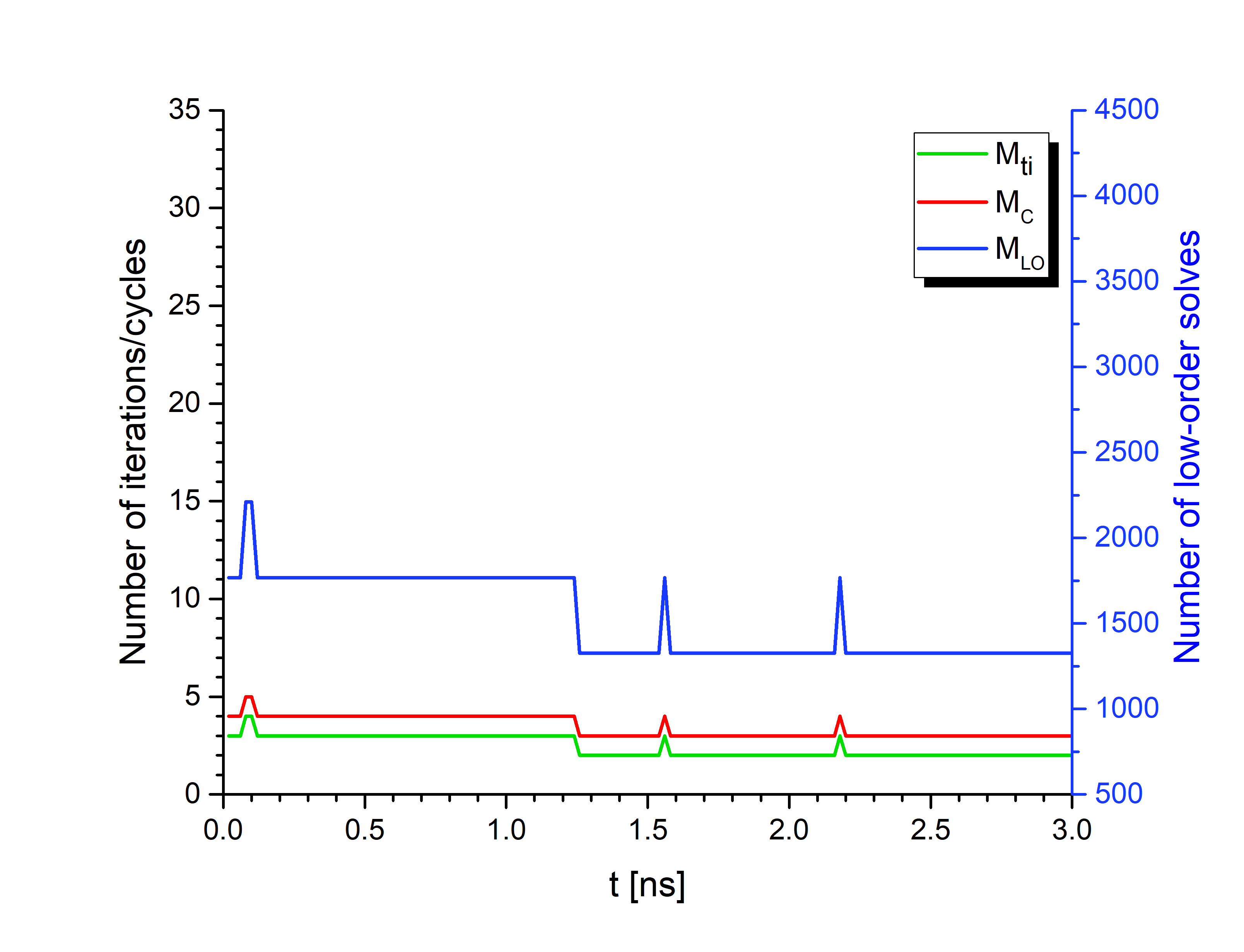}}
 \caption{\label{algorithms-256-002}  Number of transport iterations ($M_{ti}$), cycles ($M_c$), and low-order solves ($M_{lo}$)   at each time step  in the FC test with 256 groups  and $\Delta t = 2 \! \times \! 10^{-2}$ ns  over  $t\in$[0, 3\,ns].}
 \end{figure}

\begin{figure}[h!]
\centering
\subfloat[$\Gamma$=2, $n_{\nu}^{\gamma}$=256,1, $V$-cycle, $\ell_{max}$=6. \label{V-256-004-fig}]
{\includegraphics[scale=0.28]{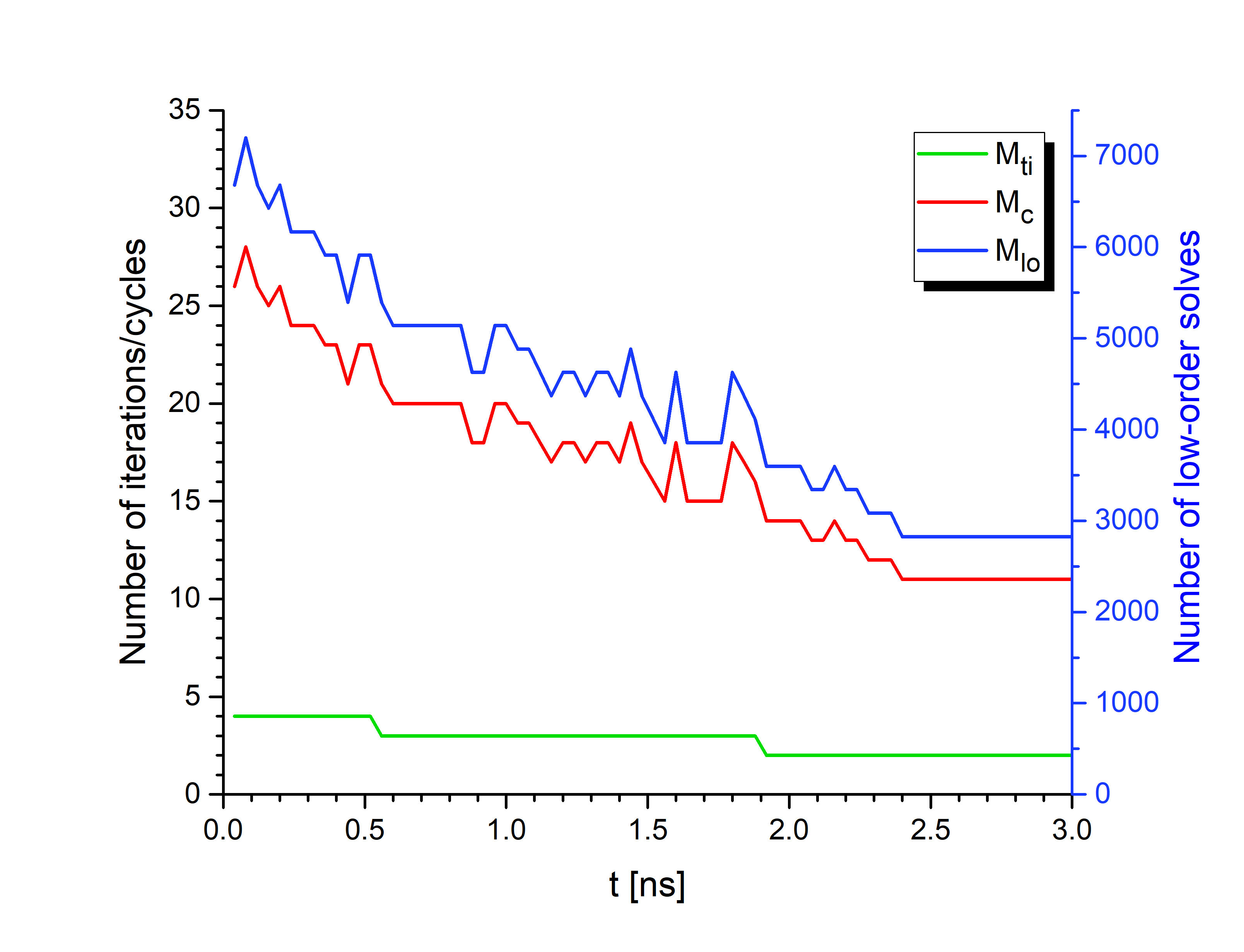}}
\hfill
\subfloat[$\Gamma$=3, $n_{\nu}^{\gamma}$=256,32,1, $W$-cycle, $\ell_{max}$=3. \label{W-256-32-004-fig}]{\includegraphics[scale=0.28]{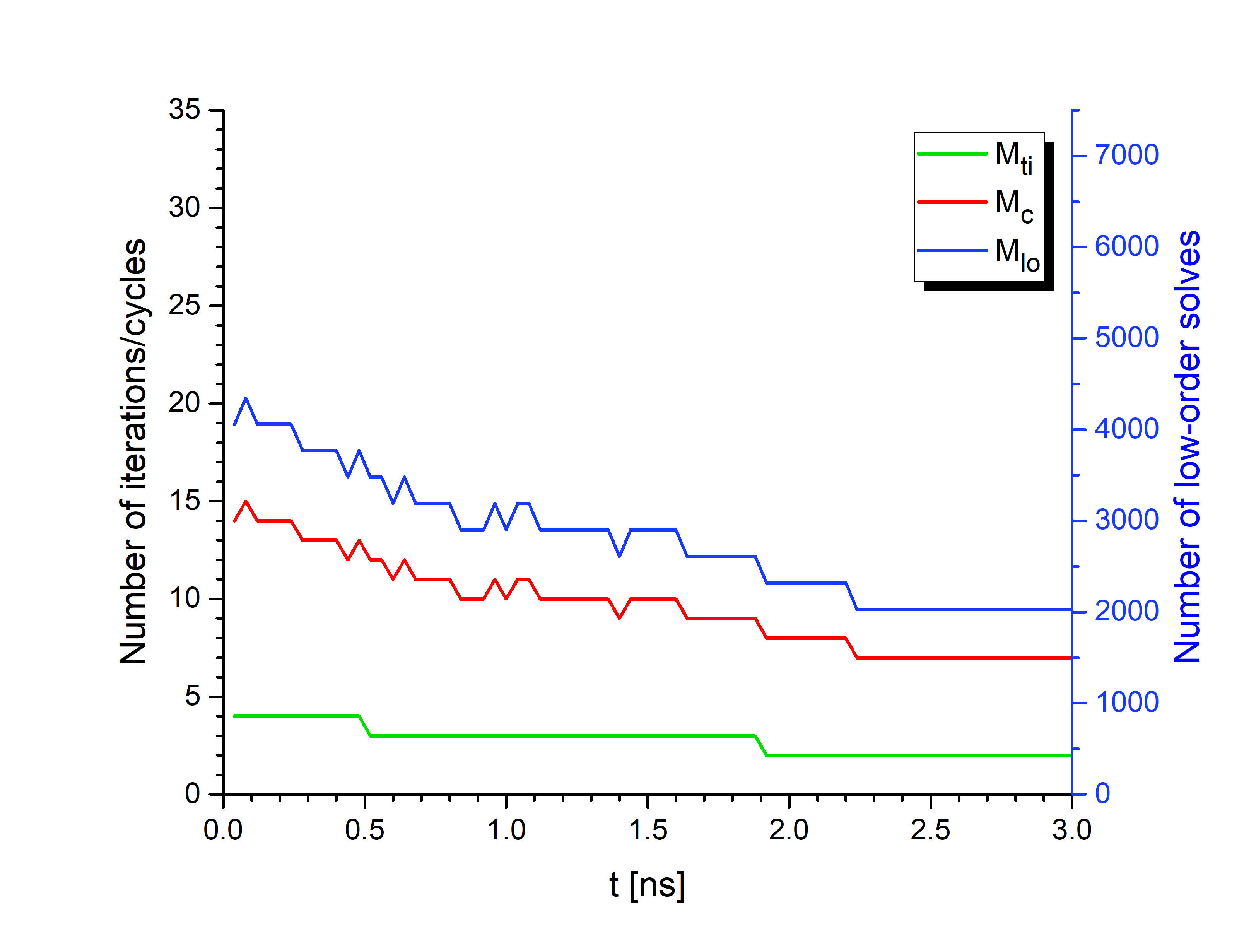}}\\
\subfloat[$F$-cycle, $\ell_{max}$=3,  $n_{\nu}^{\gamma}=256,32,16,1$, $\Gamma$=4.
\label{F-256-32-004-fig}]
{\includegraphics[scale=0.28]{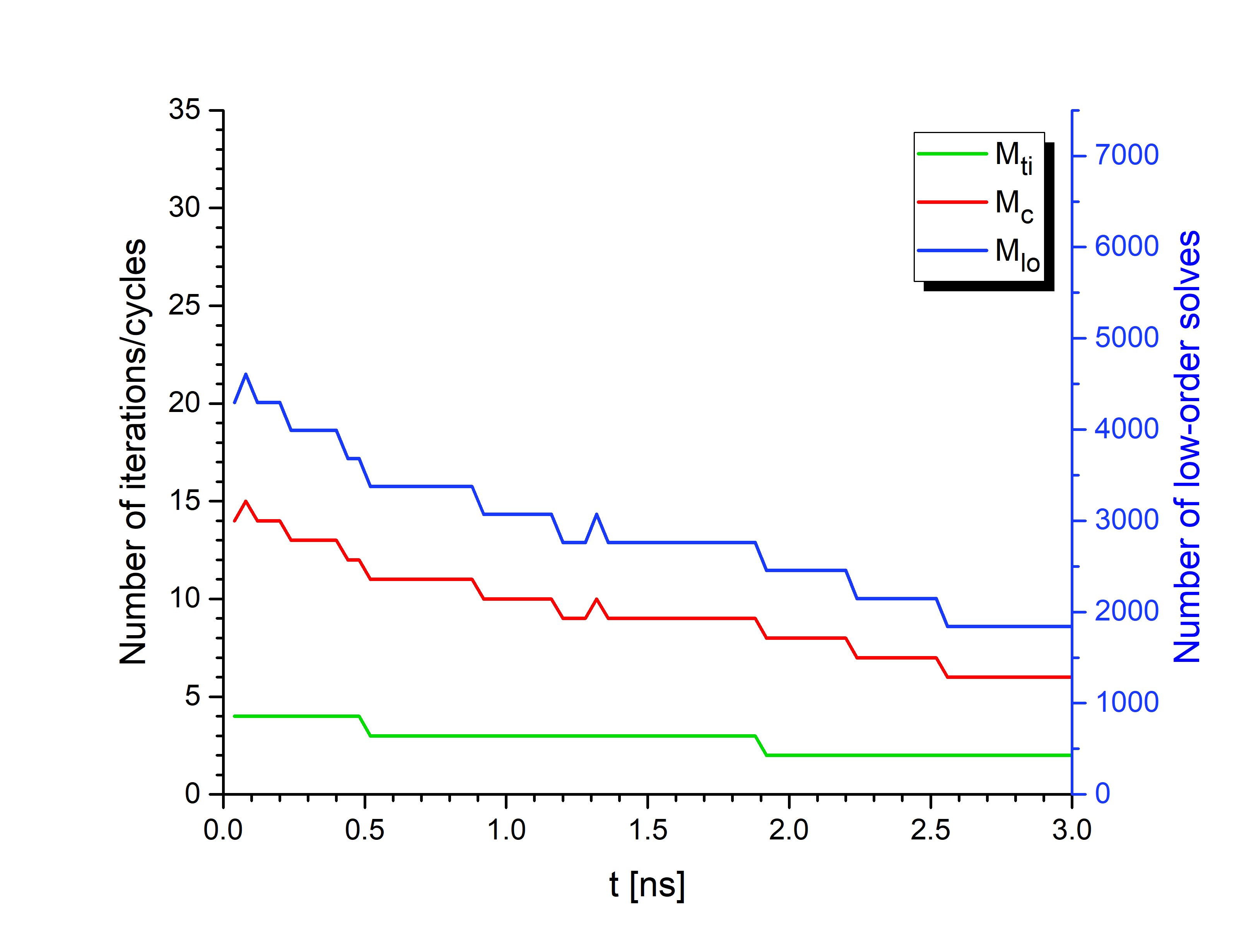}}
\hfill
\subfloat[ $\Gamma$=5, $n_{\nu}^{\gamma}$=256,32,16,4,1, $F$-cycle, $\ell_{max}$=2. \label{F-256-32-16-4-004-fig}]{\includegraphics[scale=0.28]{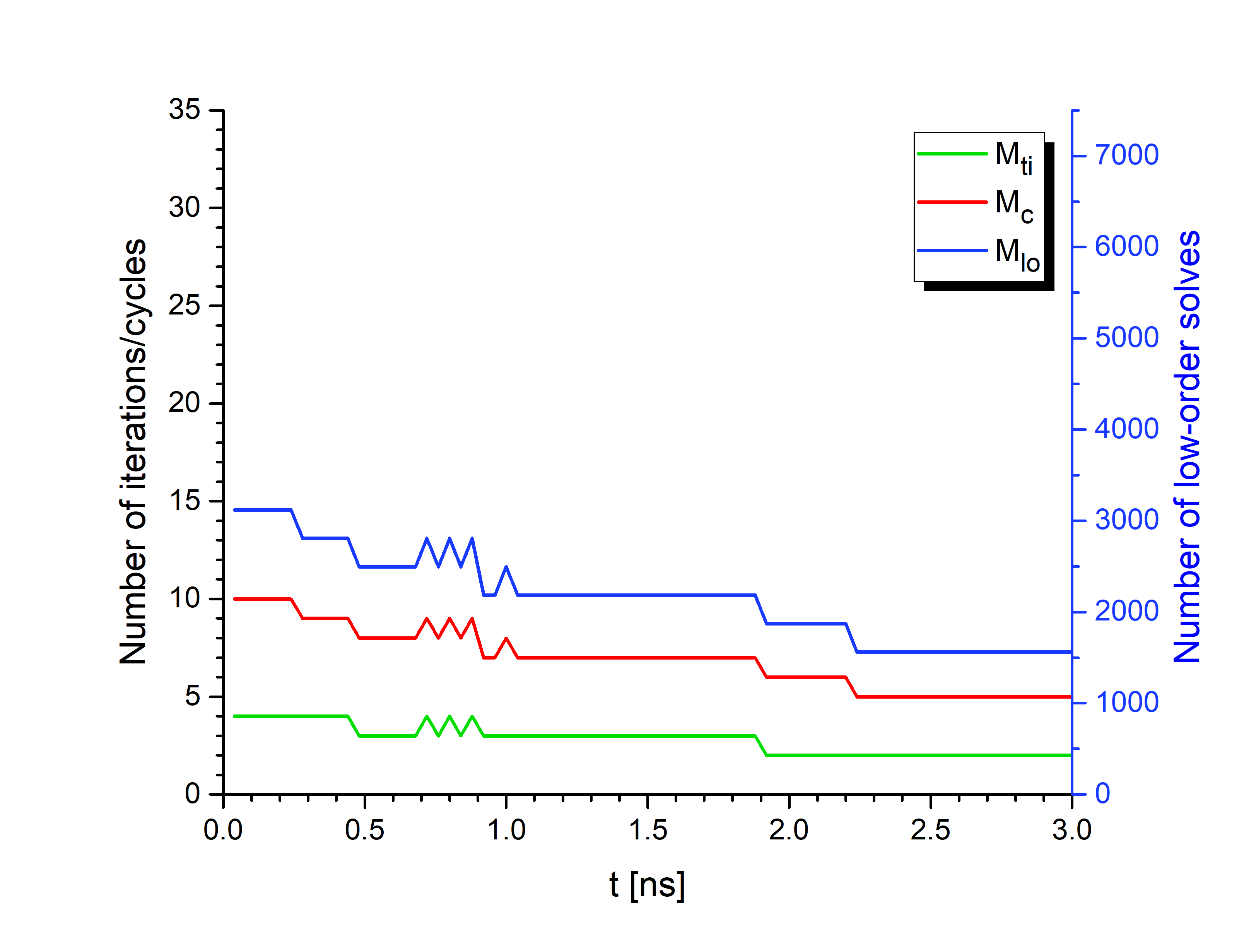}}\\
\subfloat[ $\Gamma$=6, $n_{\nu}^{\gamma}$=256,64,32,16,4,1, $F$-cycle, $\ell_{max}$=2. \label{F-256-64-32-16-4-004-fig}]{\includegraphics[scale=0.28]{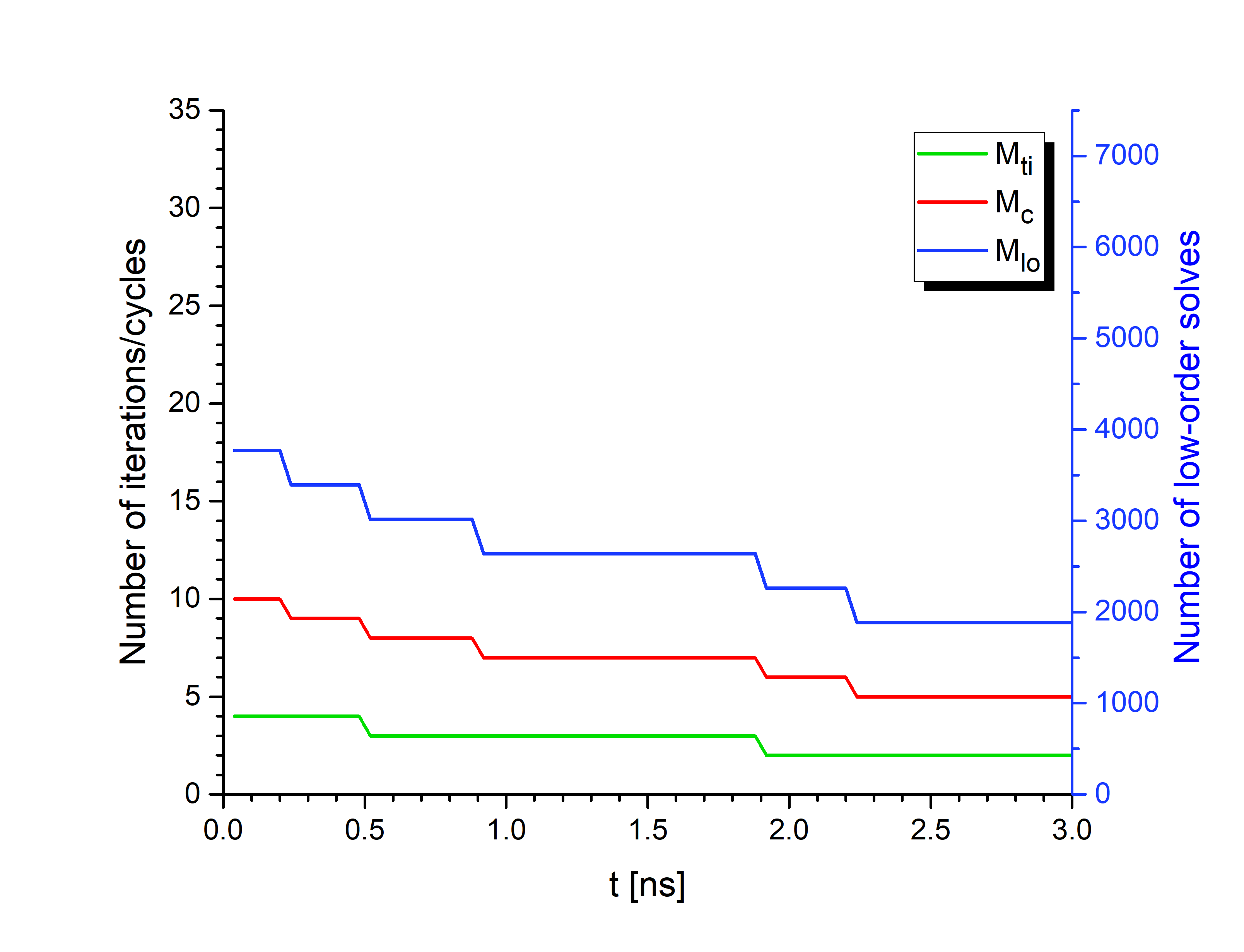}}
\hfill
\subfloat[$\Gamma$=7,\,$n_{\nu}^{\gamma}$=256,64,32,16,8,4,1,\,$F$-cycle,\,$\ell_{max}$=2. \label{F-256-64-32-16-8-4-004-fig}]{\includegraphics[scale=0.28]{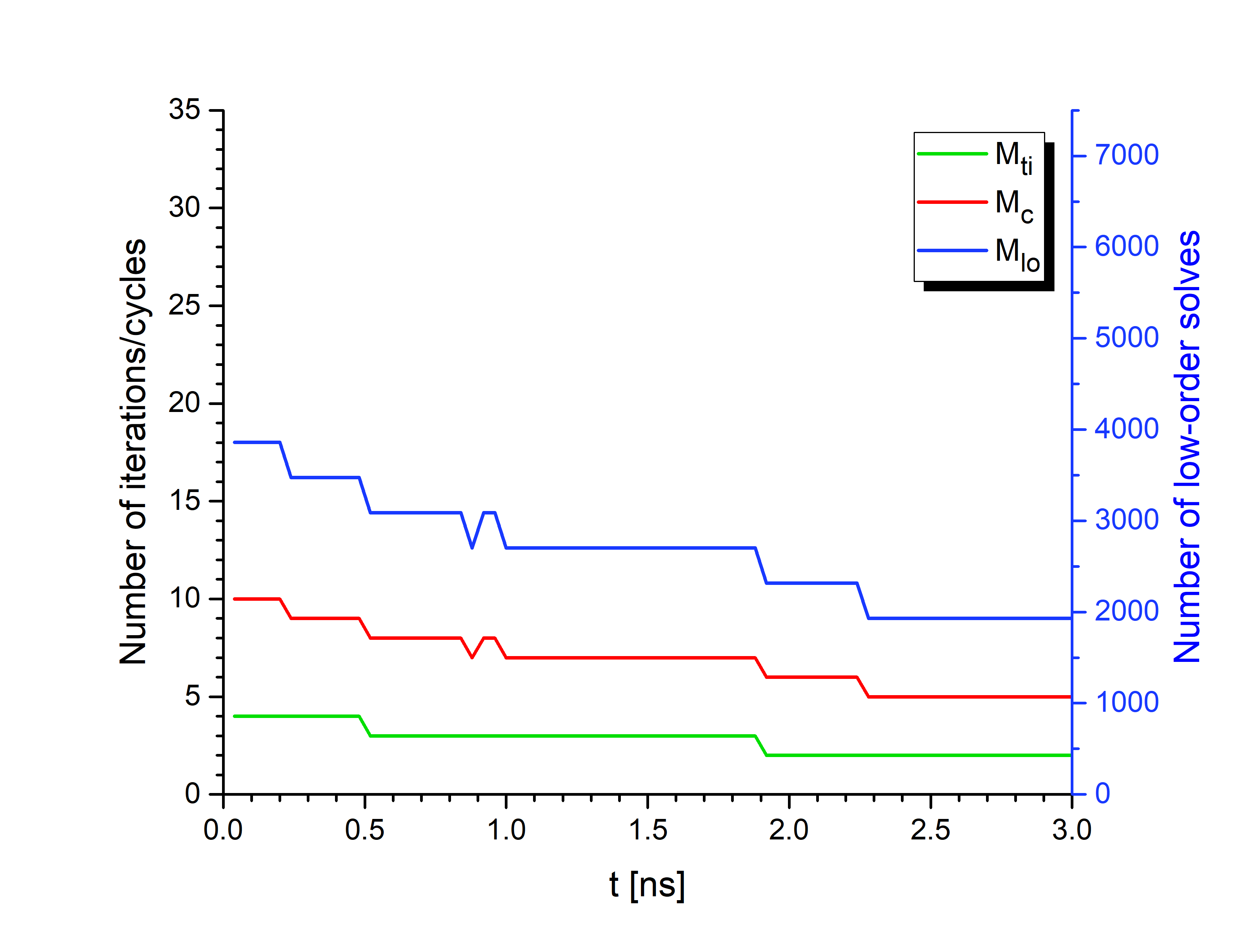}}
 \caption{\label{algorithms-256-004}  Number of transport iterations ($M_{ti}$), cycles ($M_c$), and low-order solves ($M_{lo}$)   at each time step  in the FC test with 256 groups  and $\Delta t = 4 \! \times \! 10^{-2}$ ns  over  $t\in$[0, 3\,ns].}
 \end{figure}

Figures \ref{algorithms-256-002} and  \ref{algorithms-256-004} present
 the results of the  algorithms
on selected hierarchies of grids in tests with  $\Delta t=2 \! \times \! 10^{-2}$ ns and
  $\Delta t=4 \! \times \! 10^{-2}$ ns, respectively.
These figures   show the
numbers of transport iterations ($M_{ti}$), cycles ($M_c$),
and low-order solves ($M_{lo}$)   versus time instant
for each algorithm in two test cases.
For each value of $\Gamma$,
the results are presented for such combination of the type of algorithm, $\ell_{max}$ and  set of grids that yields the smallest total number of transport iterations and total number of cycles.
The hierarchies of grids and the value of $\ell_{max}$ are indicated in the figures.
Tables  \ref{tbl-256-002} and  \ref{tbl-256-004} show  the total numbers of transport iterations ($N_{ti}$),
   cycles ($N_c$), low-order solves ($N_{lo}$)
 in tests with  $\Delta t=2 \! \times \! 10^{-2}$ ns and  $\Delta t=4 \! \times \! 10^{-2}$ ns, respectively.

Each time step is a different case for iteration algorithms.
However, there are similarities between instants over characteristic stages of the TRT problem, such as
(a) initial radiation penetration in the domain,
(b) radiation wave formation,
(c)  well-developed  wave,
and (d) approaching  steady-state regime.
The differences in evolution of the solution during these stages affect the numbers of iterations and cycles.
 These effects can be seen in Figures \ref{algorithms-256-002} and  \ref{algorithms-256-004}.
The bigger the time step the larger the change in the solution over the time  interval.
 This leads to some increase in numbers of   iterations.

The obtained results show that the   algorithms with multiple grids in frequency  significantly  reduce the number of cycles.
 We note that each algorithm has different computational costs.
They depend on (i) the number of transport iterations and (ii) the number of low-order solves, i.e. the number of times the group-wise LOQD equations are solved.
Tables  \ref{tbl-256-002} and  \ref{tbl-256-004} also show the ratio between the number of low-order solves required in the tests  by the algorithm  ($N_{lo}$) and this
number in the case of the $V$-cycle ($\Gamma=2$).

 The number of transport iterations ($N_{ti}$)  is an important factor of algorithm efficiency,
  because   their computational costs   are directly proportional to the
 number of angular directions and as well as number of frequency groups in the problem.
There is a small variation in $N_{ti}$ among the presented algorithms.
In  the test with $\Delta t=2 \! \times \! 10^{-2}$ ns,
the algorithm with the $V$-cycle on 2 grids needs 365 transport iterations.
This can be considered as a target value of $N_{ti}$ in this test.
Most of the algorithms execute  366-367 transport iterations.
The  algorithm   with  the $W$-cycle  and $\ell_{max}=2$
on the grids with $n_{\nu}^{\gamma}$=256, 32, 1
requires   362   transport iterations.
The algorithm with  the $F$-cycle  on  7 grids  with $N_{ti}$=366 executes the smallest number of low-order solves.
In  the case of   $\Delta t=4 \! \times \! 10^{-2}$ ns, the algorithm with the $V$-cycle   executes 210 transport iterations.
Almost all multigrid algorithms require  209 transport iterations in this case.
Among them the algorithm with the $F$-cycle on 6 grids has the smallest number of low-order solves.
The algorithm with the $F$-cycle on the grids with $\Gamma$=5    executes even less number of low-order solves $N_{lo}$, but needs 2 more transport iterations.

\begin{figure}[h!]
\centering
\subfloat[$\Delta t = 2 \! \times \! 10^{-2}$ ns \label{delta Ts-0.02}]{\includegraphics[scale=0.3]{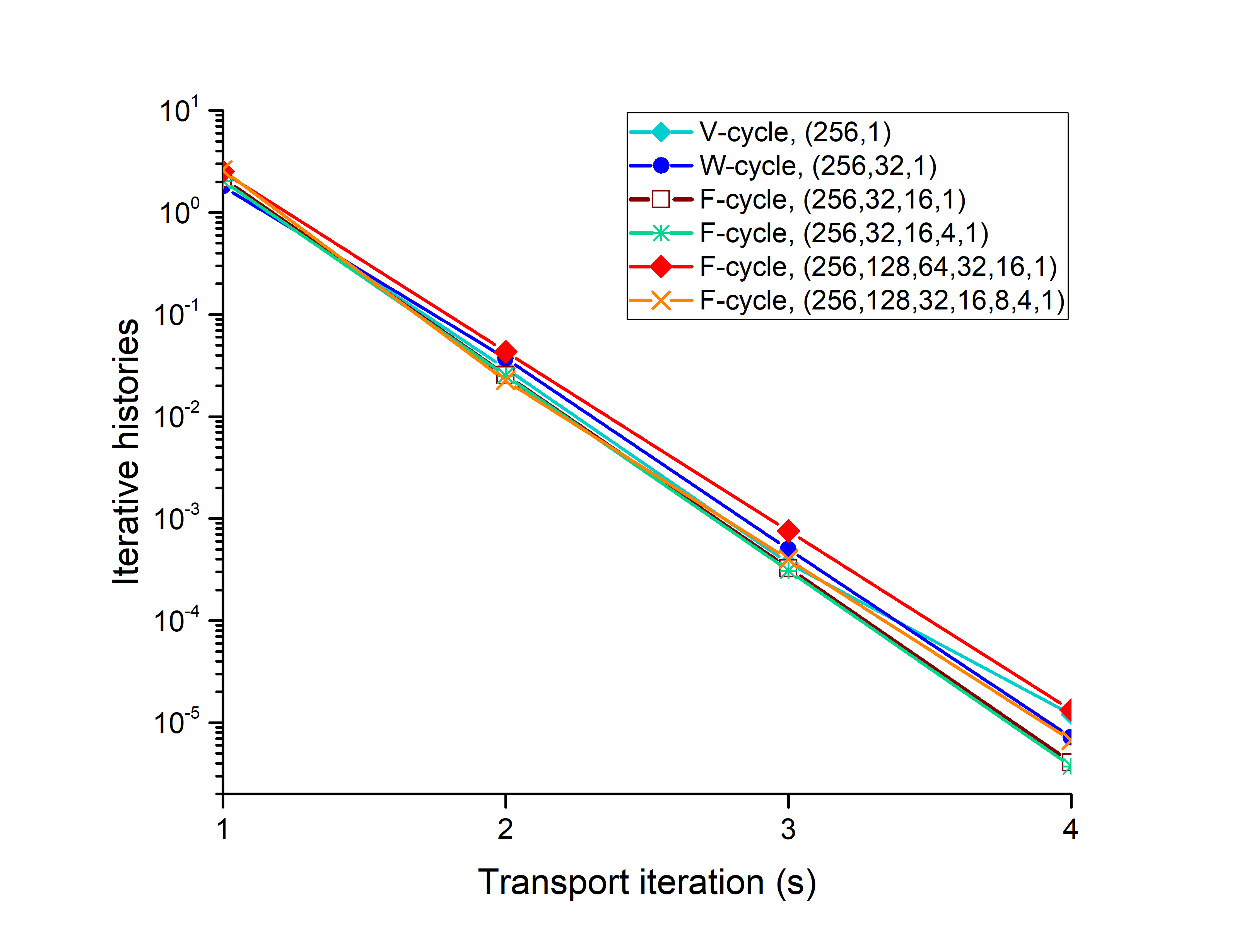}}
\subfloat[$\Delta t = 4 \! \times \! 10^{-2}$ ns \label{delta Ts-0.04}]{\includegraphics[scale=0.3]{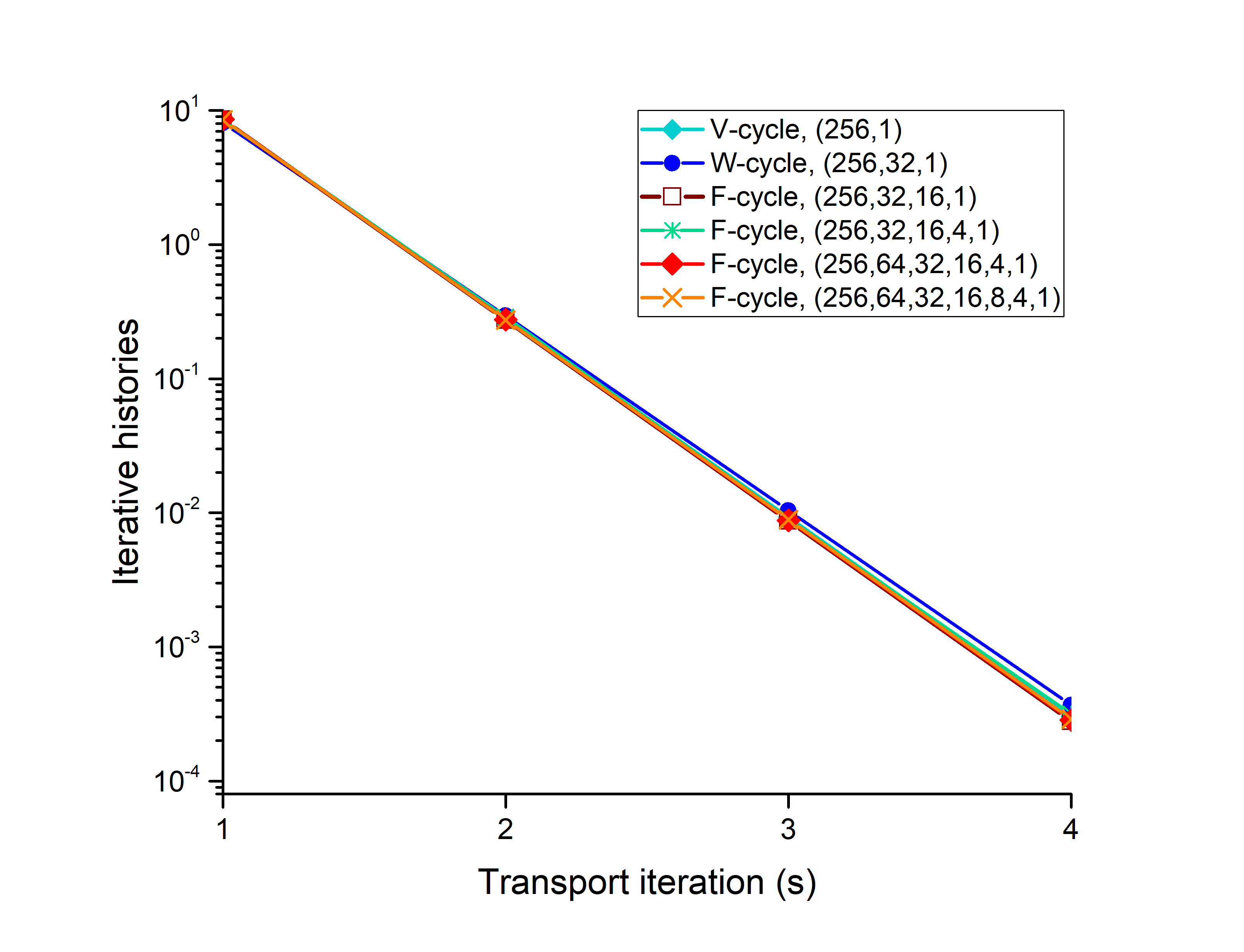}}
\caption{\label{trans-iter-hist-256-t=0.08}
Convergence of  temperature ($||\Delta T^{(s)}||_{\infty}$ [eV])  over transport iterations   at $t=8 \! \times \! 10^{-2}$ ns.}
\end{figure}

\begin{figure}[h!]
\centering
\subfloat[$\Delta t = 2 \! \times \! 10^{-2}$ ns   \label{delta Tsl-0.02}]{\includegraphics[scale=0.3]{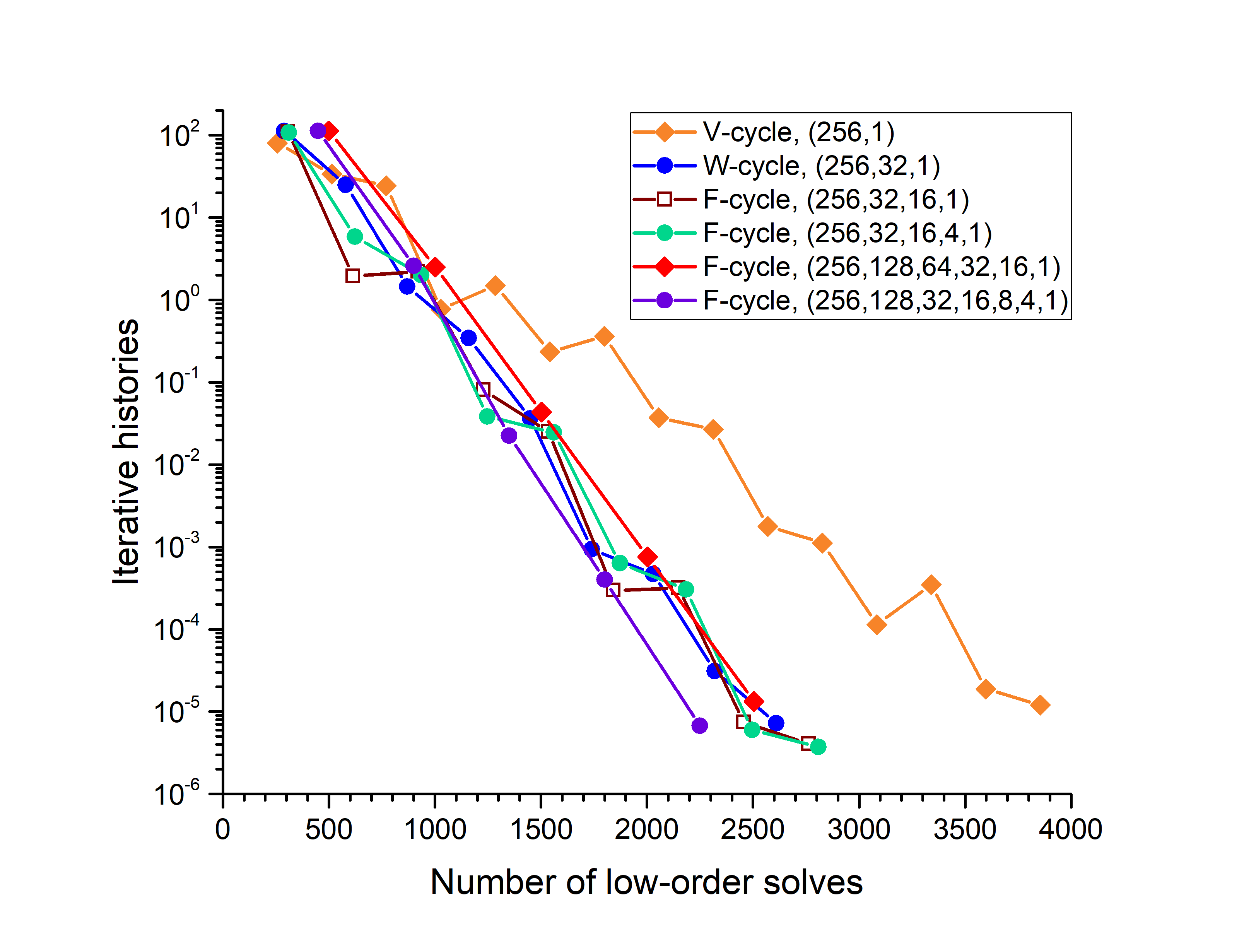}}
\subfloat[$\Delta t = 4 \! \times \! 10^{-2}$ ns \label{delta Tsl-0.04}]{\includegraphics[scale=0.3]{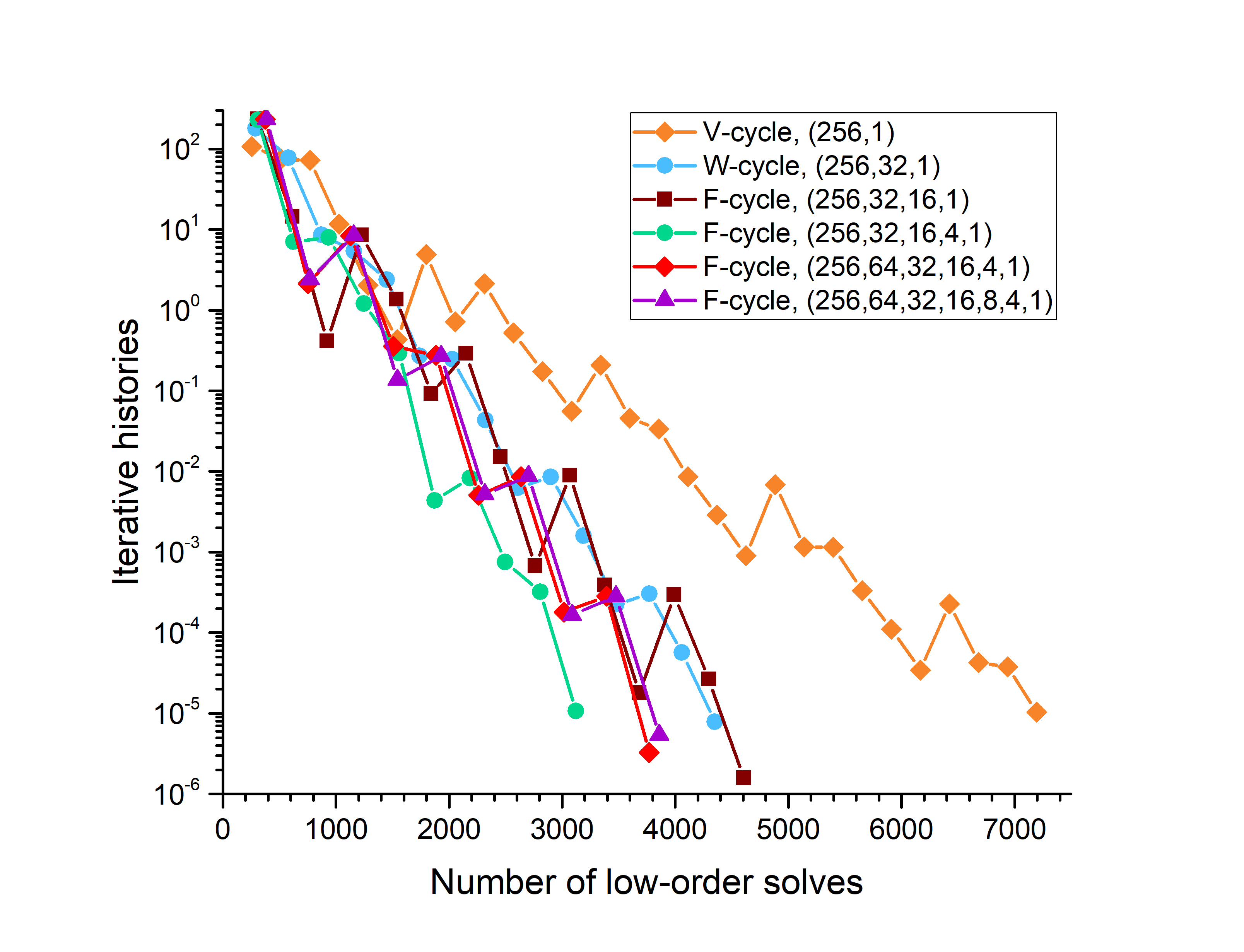}}
\caption{\label{iter-hist-256-t=0.08}
 Iterative convergence  of   temperature  ($||\Delta \tilde T^{[k]}||_{\infty}$ [eV])  over cycles    at $t=8 \! \times \! 10^{-2}$ ns.}
\end{figure}

To demonstrate convergence behavior of  different algorithms,
we use  one of initial instants, namely,  $t=8 \times 10^{-2}$ ns.
The evolution of both temperature and radiation waves is  fast  at  this early stage.
Figure \ref{trans-iter-hist-256-t=0.08} presents
convergence of  temperature   with transport iterations  ($||\Delta T^{(s)}||_{\infty}$)
in the case of $\Delta t = 2 \! \times \! 10^{-2}$ ns and $\Delta t = 4 \! \times \! 10^{-2}$ ns.
The  algorithms converge  rapidly with slightly different rates requiring the same number of transport iterations.
Figure \ref{iter-hist-256-t=0.08} shows
iterative convergence  of  temperature over cycles ($||\Delta \tilde T^{[k]}||_{\infty}$)
versus number of low-order solves.
The patterns of  convergence of inner (low-order) iterations are different.
At this instant of time, the most efficient algorithms are
(i) the $F$-cycle with $\ell_{max}=1$ on 7 grids for $\Delta t = 2 \! \times \! 10^{-2}$ ns,
and (ii)   the $F$-cycle with $\ell_{max}=3$ on 5 grids  for $\Delta t = 4 \! \times \! 10^{-2}$ ns.
They require  the smallest number of low-order solves.
We note that it takes more cycles and  low-order solves  in the case of the larger time step.
The number of transport iterations is the same for both time steps.
Thus, the computational  effort   is shifted to the  projected subspace for $\Delta t = 4 \! \times \! 10^{-2}$ ns at this instant.

\begin{table}[h!]
	\centering
	\caption{\label{tbl-256-002} Performance of    algorithms  in the FC test with 256 groups and
 $\Delta t = 2\! \times \!  10^{-2}$ ns
over 150 time steps for $t\in$[0, 3\,ns] }
\begin{tabular}{|c|c|c|c|c|c|c|c|}
  \hline
Cycle  &  $\Gamma$   & $n_{\nu}^{\gamma}$,  $\gamma \in \mathbb{N} (\Gamma)$     & $\ell_{max}$   &   $N_{ti}$  &   $N_{c}$  &  $N_{lo}$  &
$\frac{N_{lo}}{N_{lo}(V, \,\ell_{max}=4)}$   \\ \hline
$V$     &      2             &  $256,1$                    &    4                  & 365           &  1547          &   397579  &  1 \\ \hline
$W$    &      3             &   $256,32,1$              &    2                  & 362            &  901           &   261290  &  0.66 \\ \hline
$F$     &      4             &$256,32,16,1$             &    2                  & 366            &  897           &  275379  & 0.69 \\ \hline
$F$     &      5             &$256,32,16,4,1$          &    2                  &  366           &    896         &   279552 & 0.70 \\ \hline
 $F$    &      6             &$256,128,64,32,16,1$  &    1                  &  367           &    517        &   259017  & 0.65 \\ \hline
$F$     &      7             &$256,128,32,16,8,4,1$   &    1                &  366           &    516        &    232200 & 0.58  \\ \hline
  \end{tabular}
\medskip
\medskip
	\centering
	\caption{\label{tbl-256-004} Performance of    algorithms  in the FC test with 256 groups and $\Delta t = 4 \! \times \! 10^{-2}$ ns
over 75 time steps for $t\in$[0, 3\,ns] }
\begin{tabular}{|c|c|c|c|c|c|c|c|}
  \hline
Cycle  &  $\Gamma$   & $n_{\nu}^{\gamma}$,  $\gamma \in \mathbb{N} (\Gamma)$     & $\ell_{max}$   &   $N_{ti}$  &   $N_{c}$  &  $N_{lo}$  &
$\frac{N_{lo}}{N_{lo}(V, \,\ell_{max}=6)}$   \\ \hline
$V$     &      2             &  $256,1$                    &    6                  & 210           &  1262          &   324334  &  1          \\ \hline
$W$    &      3             &   $256,32,1$              &    3                  & 209            &  722           &   209380  &  0.65     \\ \hline
$F$     &      4             &$256,32,16,1$             &    3                  & 209            &  695           &  213365  & 0.66        \\ \hline
$F$     &      5             &$256,32,16,4,1$          &    2                  &  211           &    519         &   178536 & 0.55       \\ \hline
 $F$    &      6             &$256,64,32,16,4,1$     &    2                  &  209           &    516        &   194532  & 0.60        \\ \hline
$F$     &      7             &$256,64,32,16,8,4,1$   &    2                &  209           &    518        &    199948  & 0.62       \\ \hline
  \end{tabular}
\end{table}

\clearpage

\section{ \label{sec:concl} Conclusion}

This paper presented  iteration transport methods for solving TRT problems in 1D slab geometry.
They  are derived by  nonlinear projection in angular variable and frequency
and based on the MLQD method.
We developed new nonlinear iterative projection methods for TRT
that define the low-order  equations on multiple grids in photon frequency.
On each transport iteration, the multigroup LOQD equations coupled with the MEB equation
are solved iteratively by multigrid-in-frequency algorithms using the $W$- and $F$-cycles.
 The obtained results show that the new algorithms  accelerate convergence of iterations and reduce computational costs.
 The behavior of the iteration algorithms vary depending on the stage in  evolution of  temperature and radiation energy waves.  The algorithms  with different cycles and hierarchies of frequency grids  can be applied  depending on the stage of TRT  phenomenon and the value of the time step to further improve effectiveness of iteration methods.
 This kind of algorithms can be developed for TRT problems with scattering using advanced prolongation operators. Other existing iteration methods for TRT problems can take advantage in applying multigrid in frequency to develop advanced iteration techniques for problems with very large number of frequency groups.

\bibliography{dya-arxiv-2020}
\bibliographystyle{elsarticle-num}

\end{document}